\documentclass[11pt,twoside]{article}

\usepackage[latin1]{inputenc}

\usepackage[english]{babel}
\usepackage{amsfonts,amssymb,amsmath}
\usepackage{mathrsfs}
\usepackage{mathtext}
\usepackage{textcomp}
\usepackage{graphicx}
\usepackage{epsfig}
\usepackage{rotating}
\usepackage{tikz}
\usepackage{fancybox}

\usetikzlibrary{arrows,decorations.pathmorphing,backgrounds,fit}

\usepackage[toc,page]{appendix}
\usepackage[dvips, colorlinks=true,pdfstartview=FitV, linkcolor=blue,citecolor=red,urlcolor=blue]{hyperref}

\newcommand{\lra}{\longrightarrow}
\newcommand{\setm}{\setminus}
\newcommand{\ra}{\rightarrow}
\newcommand{\DS}{\displaystyle}

\newcommand{\C}{\mathbb{C}}
\newcommand{\R}{\mathbb{R}}
\newcommand{\Z}{\mathbb{Z}}

\newcommand{\N}{\mathbb{N}}
\newcommand{\E}{\mathbb{E}}
\renewcommand{\S}{\mathbb{S}}
\renewcommand{\P}{\mathrm{P}}
\newcommand{\T}{\mathrm{T}}
\newcommand{\hT}{\hat{\mathrm{T}}}

\newcommand{\A}{\mathcal{A}}
\newcommand{\Ah}{\hat{\mathcal{A}}}

\newcommand{\Ad}{\mathcal{AD}}
\newcommand{\Tr}{\mathrm{Tr}}
\newcommand{\I}{\mathcal{I}}
\newcommand{\lbd}{\lambda}

\newcommand{\Sig}{\Sigma}
\newcommand{\hSig}{\hat{\Sigma}}
\newcommand{\di}{\mathbf{d}}

\newcommand{\inter}{\mathrm{int}}
\newcommand{\leng}{\mathbf{leng}}

\newcommand{\card}{\mathrm{Card}}
\newcommand{\diam}{\mathbf{diam}}
\newcommand{\ver}{\mathrm{Ver}}

\newcommand{\Ra}{\mathrm{R}}

\newcommand{\area}{\mathbf{Area}}
\newcommand{\rk}{\mathrm{rk}}
\newcommand{\et}{\mathrm{et}}

\newcommand{\ie}{\textit{i.e. }}
\renewcommand{\geq}{\geqslant}
\renewcommand{\leq}{\leqslant}
\newcommand{\vide}{\varnothing}

\newcommand{\Hg}{\mathcal{H}(k_1,\dots,k_n)}
\newcommand{\Hgi}{\mathcal{H}_1(k_1,\dots,k_n)}

\newcommand{\transpose}[1]{{\vphantom{#1}}^{\mathit t}{#1}}

\newcommand{\MET}{\mathcal{M}^{\mathrm{et}}(\Ah,\underline{\alpha})}

\newcommand{\METTor}{\mathcal{M}^{\mathrm{et}}(\I,(2\pi,2\pi))}

\newcommand{\MSS}{\mathcal{M}(\S^2,\underline{\alpha})}
\newcommand{\MSSi}{\mathcal{M}_1(\S^2,\underline{\alpha})}
\newcommand{\MSSv}{\mathcal{M}(\S^2,\underline{\alpha})^*}
\newcommand{\MSSvi}{\mathcal{M}_1(\S^2,\underline{\alpha})^*}

\newcommand{\ST}{\mathbf{S}_\mathrm{T}}

\newcommand{\U}{\mathcal{U}}
\newcommand{\V}{\mathcal{V}}

\newcommand{\AT}{\mathbf{A}_{\mathrm{T}}}
\newcommand{\ATv}{\mathbf{A}_\mathrm{T}^*}

\newcommand{\As}{\mathbf{A}_s}
\newcommand{\Asv}{\mathbf{A}^*_s}

\newtheorem{theorem}{Theorem}[section]
\newtheorem{definition}[theorem]{Definition}
\newtheorem{proposition}[theorem]{Proposition}
\newtheorem{prop-def}[theorem]{Proposition-definition}
\newtheorem{corollary}[theorem]{Corollary}
\newtheorem{lemma}[theorem]{Lemma}

\newcommand{\dem}{\noindent {\textbf{Proof:} }}
\newcommand{\rem}{ \noindent \textbf{Remark: }}

\newcommand{\carre}{\hfill $\Box$\\}

\begin{document}
\title{\huge ENERGY FUNCTIONS ON MODULI SPACES OF FLAT SURFACES WITH ERASING FOREST}

\author{ {\large DUC-MANH NGUYEN} \\ Max-Planck-Institut für Mathematik \\ Vivatsgasse 7\\ D-53111 Bonn\\ Germany}

\date{ January 3, 2010}

\maketitle

\begin{abstract}

This paper follows on from \cite{Ng1}, in which we study flat surfaces with erasing forest, these surfaces are obtained by deforming the metric
structure of translation surfaces, and their moduli space can be viewed as some deformations of the moduli space of translation surfaces. We
showed that the moduli spaces of such surfaces are complex orbifolds, and admit a natural volume form $\mu_\Tr$. The aim of this paper is to
show that the volume of those moduli spaces with respect to $\mu_\Tr$ normalized by some energy function involving the area, and the total
length of the erasing forest, is finite. Since translation surfaces, and flat surfaces of genus zero can be viewed as special cases of flat
surfaces with erasing forest, and on their moduli space, the volume form $\mu_\Tr$ equals the usual ones up to a multiplicative constant, this
result allows us to recover some classical results of Masur-Veech, and of Thurston concerning the finiteness of the volume of the moduli space
of translation sufaces, and of the moduli space of polyhedral flat surfaces.\\

\end{abstract}

\section{Introduction}
In \cite{Ng1}, we have introduced the notion of {\em flat surface with erasing forest}. An {\em erasing forest} $\hat{A}$ in a flat surface with
conical singularities $\Sig$ is a union of disjoint geodesic trees such that

\begin{itemize}
\item[$\bullet$] the vertex set of $\hat{A}$ contains all the singularities of $\Sig$,

\item[$\bullet$] the holonomy of any closed curve which does not intersect the forest $\hat{A}$ is a translation of $\R^2$.

\end{itemize}

\noindent Note that a `generic' flat surface does not admit any erasing forest.\\

Recall that a {\em translation surface} is a flat surface with conical singularities verifying the following property: the holonomy of any
closed curve (which does not contain any singularity) is a translation. Given a translation surface $\Sig$, we can construct a flat surface with
erasing forest by deforming its metric structure as follows: first, cut off a small disk about a singular point of $\Sig$, note that by the
definition of translation surface, the cone angle at any singular point of $\Sig$ must belong to $2\pi\N$. We can modify the metric structure
inside the small disk to get a flat disk with several singular points, whose cone angles can be chosen arbitrarily, while the boundary stays
unchanged. We can then glue the disk back to $\Sig$. If the boundary is convex, then it is not hard to show that there exists a geodesic tree
inside the disk whose vertex set is the set of singularities. Carrying out this operation for all the singular points of $\Sig$, we get a new
flat surface $\Sig'$ together with a family of geodesic trees. By construction, the union of these trees is an erasing forest of $\Sig'$.\\

\noindent A translation surface is a particular flat surface with erasing forest, where each tree in the erasing forest is just a singular
point. A flat surface of genus zero can also be viewed as a flat surface with erasing forest, since there always exists a geodesic tree on this
surface whose vertex set is the set of singularities, and the complement of such a tree is just a topological disk.\\

Given a flat surface $\Sig$ with an erasing forest $\hat{A}$, a {\em parallel vector field} on $\Sig$ is a vector field defined on the
complement of the erasing forest $\hat{A}$ which is invariant by the parallel transport. In a local chart of the flat metric structure, the
integral lines of such a field are parallel.  On any (connected) flat surface with erasing forest, such vector fields always exist, they are
uniquely determined by a tangent vector at a fixed point in the complement of the erasing forest.\\

\medskip

Given an integer $g\geq 0$, and positive real numbers $\alpha_1,\dots,\alpha_n,$ verifying

$$\sum_{i=1}^n\alpha_i=(2g+n-2)2\pi,$$

\noindent let us fix a family $\Ah=\{\A_1,\dots,\A_m\}$ of topological trees such that the total number of vertices of the trees in $\Ah$ is
$n$, and choose a numbering on the set of vertices of $\Ah$. Note that we consider an isolated point as a special tree. Let $\underline{\alpha}$
denote the vector $(\alpha_1,\dots,\alpha_n)$, and $\MET$ denote the set of triples $(\Sig,\hat{A}, \xi)$ where

\begin{itemize}

\item[$\bullet$] $\Sig$ is a closed, connected, oriented flat surface of genus $g$ with cone singularities,

\item[$\bullet$] $\hat{A}$ is an erasing forest in $\Sig$ consisting of $m$ geodesic trees $A_1,\dots,A_m$, we also suppose that the trees and
vertices of $\hat{A}$ are numbered so that $A_j$ is isomorphic to $\A_j$ (as topological trees), and by those isomorphisms, the $i$-th vertex of
$\Ah$ is mapped to a singular point with cone angle $\alpha_i$.

\item[$\bullet$] $\xi$ is a unitary parallel vector field defined on $\Sig\setm\hat{A}$.

\end{itemize}

\noindent In \cite{Ng1}, we proved that $\MET$ has a structure of analytic complex orbifold of dimension

$$\left\{%
\begin{array}{ll}
    2g+n-1 & \hbox{ if $\alpha_i\in 2\pi\N,\; \forall i=1,\dots,n$, } \\
    2g+n-2 & \hbox{ otherwise,} \\
\end{array}%
\right.$$

\noindent together with a natural volume form $\mu_{\Tr}$. Remark that, as all the trees in the erasing forest shrink to points, a flat surface
with erasing forest becomes a translation surface. Therefore, $\MET$ can be viewed as a deformation of some stratum of the moduli space of
translation surfaces.\\

\medskip

Consider the following function on $\MET$

$$\begin{array}{crcc}
  \mathcal{F}^\mathrm{et}: & \MET  & \lra & \R \\
   & (\Sig,\hat{A},\xi) & \longmapsto & \exp(-\mathbf{Area}(\Sig)-\ell^2(\hat{A})) \\
\end{array}$$

\noindent where $\ell(\hat{A})$ is the total length of the trees in $\hat{A}$. In what follows, we will call a topological tree which is not a
point a {\em non-trivial tree}. The main result of this paper is the following

\begin{theorem}\label{FVTh}

If at least one of the trees in the family $\Ah$ is non-trivial, then the integral of the function $\mathcal{F}^{\mathrm{et}}$ over $\MET$ with
respect to $\mu_\Tr$ is finite:

\begin{equation}\label{FinIneq}
\int_{\MET} \mathcal{F}^\mathrm{et}d\mu_\Tr<\infty.
\end{equation}

\end{theorem}

\rem

\begin{itemize}

\item[$\bullet$] The integral (\ref{FinIneq}) is still finite if we multiply the total length of the erasing forest by a parameter $\epsilon
>0$, that is the statement of Theorem \ref{FVTh} is also true for functions $\mathcal{F}^\mathrm{et}_\epsilon : (\Sig,\hat{A},\xi) \mapsto
\exp(-\mathbf{Area}(\Sig)-\epsilon\ell^2(\hat{A}))$, with $\epsilon >0$.

\item[$\bullet$] Let $e_1,\dots,e_{n-m}$ denote the edges of the trees in the forest $\hat{A}$,  and $\ell(e_i)$ denote the length of $e_i$,
then the integral (\ref{FinIneq}) is also finite if we replace $\mathcal{F}^\mathrm{et}$ by the function

$$\widetilde{\mathcal{F}}^\mathrm{et} : (\Sig,\hat{A},\xi) \mapsto \exp(-\mathbf{Area}(\Sig)-\sum_{i=1}^{n-m} \ell^2(e_i)).$$

The proofs for $\mathcal{F}^\mathrm{et}_\epsilon$ and $\widetilde{\mathcal{F}}^\mathrm{et}$ are the same as the proof for
$\mathcal{F}^\mathrm{et}$.

\end{itemize}

In the case where all the trees in $\Ah$ are points, the space $\MET$ is identified to a stratum $\Hg$ of the moduli space of Abelian
differentials on Riemann surfaces of genus $g$, and we have

$$\mathcal{F}^\mathrm{et}(\Sig,\hat{A},\xi)=\exp(-\mathbf{Area}(\Sig)).$$

\noindent The similar result for this case can be proved as a consequence of Theorem \ref{FVTh}, that is

\begin{theorem}\label{FVprA}
We have
\begin{equation}\label{FinIneqA}
\int_{\Hg} \exp(-\mathbf{Area}(.)) d\mu_\Tr< \infty.
\end{equation}

\end{theorem}

\noindent Note that the assumption that at least one of the trees in the forest is not a point is crucial for the proof of Theorem \ref{FVTh},
hence, Theorem \ref{FVprA} cannot be considered as a particular case of Theorem \ref{FVTh}.\\

Let $\Hgi$ denote the subset of $\Hg$ consisting of surfaces of unit area. Let $\mu^1_\Tr$ denote the volume form on $\Hgi$ which is induced by
$\mu_\Tr$. A direct consequence of Theorem \ref{FVprA} is the following

\begin{corollary}\label{FVcorA}
The total measure $\DS{\mu^1_\Tr(\Hgi)}$ is finite.
\end{corollary}

\dem Identifying $\Hg$ to $\Hgi\times\R^*_+$, and we can write $\DS{d\mu_\Tr=t^sd\mu^1_\Tr dt}$, where  $s=\dim_\R \Hgi$ which is odd.
Therefore, we have

\begin{eqnarray*}
\int_{\Hg}\exp(-\mathbf{Area}(.))d\mu_\Tr &= & \int_{\Hgi}\int_{0}^{+\infty}t^se^{-t^2}dtd\mu^1_\Tr,\\
        &= & \frac{1}{2}(\frac{s-1}{2})!\int_{\Hgi}d\mu^1_\Tr\\
\end{eqnarray*}

\noindent and the corollary follows. \carre

On the space $\Hg$, we have (see \cite{MasTab}, \cite{Zor}) a ``natural" volume form $\mu_0$ which is defined by the period mapping, let
$\mu^1_0$ denote the volume form on $\Hgi$ which is induced by $\mu_0$. In \cite{Ng1}, we proved that $\DS{\mu_\Tr=\lambda\mu_0}$, where
$\lambda$ is a constant on each connected component of $\Hg$. By a well known result of Kontsevich-Zorich \cite{KonZo}, we know that $\Hg$ has
at most three connected components, thus Corollary \ref{FVcorA} is equivalent to the classical result of Masur-Veech  stating that the volume of
$\Hgi$ with respect to $\mu^1_0$ is finite.\\

Let us now consider flat surfaces of genus zero, that is  flat surfaces homeomorphic to the sphere $\S^2$. Fix  $n, \; n\geq 3$, positive real
numbers $\alpha_1,\dots,\alpha_n$ verifying

$$\sum_{i=1}^n \alpha_i=(n-2)2\pi.$$

\noindent Let $\underline{\alpha}$ denote the $n$-uple $(\alpha_1,\dots,\alpha_n)$, and  $\MSSv$ denote the moduli space of flat surfaces of
genus zero having exactly $n$ singular points with cone angles $\alpha_1,\dots,\alpha_n$. Let $\MSS$ denote the product space $\MSSv\times
\S^1$. \\

\noindent Given a point $(\Sig,e^{\imath\theta})$ in $\MSS$, it is not difficult to see that there always exists an erasing forest consisting of
only one geodesic tree $A$ in $\Sig$, therefore, a neighborhood of $(\Sig,e^{\imath\theta})$ in $\MSS$ can be identified to an open set in
$\MET$, where the family $\Ah$ contains only one tree which is isomorphic to $A$. We also get a volume form $\mu_{\Ah}$ on a neighborhood of
$(\Sig,e^{\imath\theta})$ which, {\em a priori}, depends on a choice of the erasing tree $A$. In \cite{Ng1}, we showed that the volume form
$\mu_{\Ah}$ actually does not depend on the choice of the tree $A$, therefore, we get a well defined volume form $\mu_\Tr$ on $\MSS$. Using
Theorem \ref{FVTh}, we will prove

\begin{theorem}\label{FVprB}
The integral of the function $(\Sig,e^{\imath\theta})\longmapsto \exp(-\mathbf{Area}(\Sig))$ over $\MSS$ with respect to $\mu_\Tr$ is finite:

\begin{equation}\label{FinIneqB}
\int_{\MSS}e^{-\area}d\mu_\Tr<\infty.
\end{equation}
\end{theorem}

\noindent Let $\MSSvi$ denote the subset of $\MSSv$ consisting of surfaces of unit area. The volume form $\mu_\Tr$ on $\MSS$ induces a volume
form $\hat{\mu}^1_\Tr$ on $\MSSvi$. The same arguments as in Corollary \ref{FVcorA} show

\begin{corollary}\label{FVcorB}

The volume of $\MSSvi$ with respect to $\hat{\mu}_\Tr^1$  is finite.

\end{corollary}

In the case where $\alpha_i < 2\pi$, for $i=1,\dots,n$, Thurston \cite{Thu} showed that $\MSSvi$ can be equipped with a complex hyperbolic
metric structure with finite volume. In \cite{Ng1}, it is showed that $\hat{\mu}^1_\Tr=\lambda\mu_{\mathrm{Hyp}}$, where $\lambda$ is a
constant, and $\mu_\mathrm{Hyp}$ is the volume form induced by the complex hyperbolic metric. Therefore Theorem \ref{FVprB} can be considered as
a generalization of the Thurston's result. It is also worth noticing that a similar result to Corollary \ref{FVcorB} has been proved in
\cite{Vee93}.\\

In the next section, we recall the definitions of the local charts for $\MET$, and the construction of the volume form $\mu_\Tr$. In Section
\ref{TorSect}, we will give the proof of Theorem \ref{FVTh} in a simple case.  The proof of Theorem \ref{FVTh} for the general case will be
given in Section \ref{FVThPrf}, and subsequently the proof of Theorem \ref{FVprA}, and Theorem \ref{FVprB} will be given in Section
\ref{prAPrf}, and Section \ref{prBPrfSect}.\\

\noindent {\bf Acknowledgements:} The author would like to express his gratitude toward François Labourie for the guidance, and for the
encouraging discussions, which are indispensable for the accomplishment of this work. This manuscript is written during the author's stay at
Max-Planck-Institut für Mathematik in Bonn, the author is  thankful to the Institute for its hospitality.\\

\section{Local charts and volume form on $\MET$}

In this section, we recall the definitions of local charts, and of the volume form $\mu_\Tr$ on $\MET$ as well as $\MSS$, details of proofs are
given in \cite{Ng1}.\\

\noindent Let $(\Sig,\hat{A},\xi)$ be a point in $\MET$. A geodesic triangulation $\T$ of $\Sig$ is said to be {\em admissible} if its
$1$-skeleton contains the forest $\hat{A}$. Given such a triangulation, we construct a local chart for $\MET$ in a neighborhood of
$(\Sig,\hat{A},\xi)$ as follows: first, cut open the surface $\Sig$ along the trees of $\hat{A}$, we then get a flat surface $\hat{\Sig}$ with
piecewise geodesic boundary  together with a geodesic triangulation $\hat{\T}$.\\

\noindent We choose an orientation for every (geometric) edge in the $1$-skeleton of $\hat{\T}$. Map each triangle of $\hat{\T}$ isometrically,
and preserving the orientation into $\R^2$ such that the parallel vector field $\xi$ is identified to the constant vertical vector field $(0,1)$
of $\R^2$. We can then associate to each oriented edge $e$ in the $1$-skeleton of  $\hat{\T}$ a well-defined complex number $z(e)$. The complex
numbers associated to edges of $\hat{\T}$ are obviously related, namely

\begin{itemize}
\item[$\bullet$] If $e_i,e_j,e_k$ are the edges of $\hat{\T}$ that bound a triangle then

\begin{equation}\label{TriaEq}
\pm z(e_i)\pm z(e_j) \pm z(e_k)=0
\end{equation}

where the signs are chosen according to the orientation of $e_i,e_j$, and $e_k$.

\item[$\bullet$] If $(e,\bar{e})$ is a pair of edges in the boundary of $\hat{\Sig}$ which arise from the same edge $\tilde{e}$ of a tree in
$\hat{A}$, then

\begin{equation}\label{BdrEq}
\pm z(\bar{e})\pm e^{\imath\theta}z(e)=0
\end{equation}

where $\theta$ is the rotation angle of the holonomy of a closed curve in $\Sig$ meeting $\hat{A}$  at only one point in $\tilde{e}$
transversely, $\theta$ is determined up to sign by the angles $(\alpha_1,\dots,\alpha_n)$, and the tree that contains $\tilde{e}$.

\end{itemize}

Let $N_1$ and $N_2$ be the number of edges and the number of triangles of $\hat{\T}$ respectively. Simple computations show that

$$N_1=3(2g+m-2)+4(n-m), \text{ and } N_2=2(2g+m-2)+2(n-m).$$

\noindent The complex numbers associated to the edges of $\hat{\T}$ give us a vector $Z$ in $\C^{N_1}$. The arguments above show that the
coordinates of $Z$ satisfy a system $\ST$ of linear equations consisting of

\begin{itemize}
\item[$\bullet$] $N_2$ equations of type (\ref{TriaEq}) which will be called {\em triangle equations}, and

\item[$\bullet$] $n-m$ equations of type (\ref{BdrEq}) which will be called {\em boundary equations}

\end{itemize}


Let

$$\AT: \C^{N_1}\ra\C^{N_2+(n-m)}$$

\noindent be the complex linear map which is defined in the canonical bases of $\C^{N_1}$ and $\C^{N_2+(n-m)}$ by the matrix whose entries are
coefficients of the system $\ST$. Note that every entry of the matrix of $\AT$ is either $0$, or a complex number of module $1$. We then have a
map $\Psi_\T$ defined in a neighborhood of $(\Sig,\hat{A},\xi)$ with image in $\ker\AT$, which associates to any point $(\Sig',\hat{A}',\xi')$
close to $(\Sig,\hat{A},\xi)$ a vector in $\ker\AT$ whose coordinates arise from an admissible triangulation $\T'$ of $\Sig'$ isomorphic to
$\T$. It turns out that $\Psi_\T$ is a local chart for $\MET$, as a consequence $\dim_\C \MET=\dim_\C \ker\AT$, and we have

$$\dim_\C \MET=N_1-\mathrm{rk}(\ST)=\left\{%
\begin{array}{ll}
    2g+n-1 & \hbox{ if $\alpha_i\in 2\pi\N, \; \forall i=1,\dots,n$,} \\
    2g+n-2 & \hbox{ otherwise}. \\
\end{array}
\right.$$

Using $\AT$, we define a volume form $\nu_\T$ on $\ker\AT$ as follows:

\begin{itemize}

\item[$\bullet$] If $\dim \MET=2g+n-1$, or equivalently $\mathrm{rk}\AT=N_2+(n-m)-1$, then $\nu_\T$ is the volume form on $\ker\AT$ which is
induced by the Lebesgue measures of $\C^{N_1}, \C^{N_2+(m-n)},$ and  $\C$ via the following exact sequence

\begin{equation}\label{ExSeq1}
0\lra \ker\AT \hookrightarrow \C^{N_1}\overset{\AT}{\lra} \C^{N_2+(n-m)}\overset{\mathbf{s}}{\lra}\C\lra 0
\end{equation}

where $\mathbf{s}$ is a linear form on $\C^{N_2+(n-m)}$ of the form

$$\mathbf{s}(z_1,\dots,z_{N_2+n-m})=\pm z_1\pm\dots\pm z_{N_2+n-m}.$$

\item[$\bullet$] If $\dim \MET=2g+n-2$, or equivalently $\mathrm{rk}\AT=N_2+n-m$, then $\nu_\T$ is the volume form which is induced by the
Lebesgue measures of $\C^{N_1}$, and $\C^{N_2+n-m}$ via the exact sequence

\begin{equation}\label{ExSeq2}
0 \lra \ker\AT \hookrightarrow \C^{N_1} \overset{\AT}{\lra} \C^{N_2+n-m} \lra 0
\end{equation}

\end{itemize}

\noindent Let $\mu_\T$ denote $\Psi^*_\T\nu_\T$, then $\mu_\T$ is a volume form defined in a neighborhood of $(\Sig,\hat{A},\xi)$. It turns out
that the volume form $\mu_\T$ does not depend on the choice of the triangulation $\T$, thus we get a well defined volume form on $\MET$ which is
denoted by $\mu_\Tr$.\\

Recall that $\MSSv$ is the moduli space of flat surfaces of genus zero having exactly $n$ singularities with cone angles given by
$\underline{\alpha}=(\alpha_1,\dots,\alpha_n)$. Let $\Sig$ be a point in $\MSSv$, then there exists a geodesic tree $A$ on $\Sig$ whose vertex
set is the set of singular points, such a tree is by definition an erasing forest of $\Sig$. As a consequence, a neighborhood of a point
$(\Sig,e^{\imath\theta})$ in $\MSS=\MSSv\times\S^1$ can be identified to a neighborhood of a point $(\Sig,A,\xi)$ in $\MET$, where $\Ah$
contains only one tree, which is isomorphic to $A$. We can then use the same method as above to define local charts, and the volume form
$\mu_{\Tr}$ for $\MSS$.\\

\noindent Note that in this case there always exist indices $i \in \{1,\dots,n\}$ such that $\alpha_i \not\in 2\pi\N$, since we must have
$\DS{\alpha_1+\dots+\alpha_n=(n-2)2\pi}$. It follows that $\DS{\dim_\C \MSS=n-2}$, and $\mu_{\Tr}$ is defined by the exact sequence
(\ref{ExSeq2}). The fact that $\mu_\Tr$ is well-defined follows from the observation that any two geodesic triangulations of $\Sig$ whose vertex
sets coincide with the set of singular points of $\Sig$ can be transformed, one into the other, by a sequence of elementary moves (see
\cite{Ng1}, Definition 6.1).\\

\section{Case of flat tori with marked geodesic segments}\label{TorSect}

In this section, we prove Theorem \ref{FVTh} for the case $g=1, n=2, m=1, \alpha_1=\alpha_2=2\pi$, and $\Ah=\{\I\}$ where $\I$ is a segment. Via
this simple case, we would like to illustrate the strategy of the proof of Theorem \ref{FVTh} in the general case. An element of $\METTor$ is a
triple $(\Sig,I,\xi)$, where $\Sig$ is a flat torus (without singularity), $I$ is an oriented geodesic segment in $\Sig$ with distinct
endpoints, the orientation of $I$ arises from a numbering of its endpoints, and $\xi$ is a unitary parallel vector field on $\Sig$.  Note that

$$\dim_\C \METTor=3.$$

\noindent Given an element $(\Sig, I, \xi)$ in $\METTor$, let $p$ and $q$  denote the endpoints of $I$ so that the orientation of $I$ is from
$p$ to $q$. Let us start by showing that one can always cut the torus $\Sig$ into two cylinders such that one of which contains $I$. This will
allows us to get a domain in $\C^3$ which covers a full measure subset of $\METTor$.

\begin{lemma} \label{FVTorL1}
There always exists a pair of parallel simple closed geodesic $\gamma_p, \gamma_q$ of $\Sig$ such that

$$\gamma_p\cap I=\{p\}, \text{ and } \gamma_q\cap I=\{q\}.$$
\end{lemma}

\dem  Choose a direction $\theta$ which is not parallel to $I$, and let $(\psi^\theta_t), \; t\in \R$, denote the geodesic flow on $\Sig$ in
this direction. Observe that there exists $t >0$ such that

\begin{equation}\label{TorCd1}
 \psi^\theta_t(I)\cap I \neq \vide
\end{equation}

\noindent since otherwise, the area of the stripe swept out by $(\psi^\theta_t)_{t>0}(I)$ would tend to infinity. Let $t_0>0$ be the first time
such that (\ref{TorCd1}) holds.  By definition, there exists a closed parallelogram $P$ in $\R^2$ with two horizontal sides, and an isometric
immersion $\DS{\varphi: \mathrm{P}\lra \Sig}$, whose  restriction to $\inter(\P)$ is an embedding, which maps the lower horizontal side of $\P$
to $I$, and the upper horizontal side of $\P$ to $\psi^\theta_{t_0}(I)$. Since the segments $I$ and $\psi_{t_0}^\theta(I)$ are parallel, and
have the same length, their intersection contains at least one endpoint of $I$. Without loss of generality, we can assume that

$$p\in I\cap\psi^\theta_{t_0}(I).$$

\noindent Consequently, $\varphi^{-1}(p)$ contains exactly two points, one in lower horizontal side, and the other in the upper horizontal side
of $\P$.\\

\noindent Let $s$ be the geodesic segment in $\P$ joining two points in $\varphi^{-1}(p)$, then $\gamma_p=\varphi(s)$ is a closed geodesic in
$\Sig$ which intersects $I$ only at $p$. The closed geodesics parallel to $\gamma_p$ which intersect $I$ fill out a cylinder whose boundary
consists of $\gamma_p$, and the closed geodesic parallel to $\gamma_p$ passing through $q$, we denote this geodesic by $\gamma_q$. By
construction, $\gamma_p$ and $\gamma_q$ satisfy the required condition of the lemma. \carre

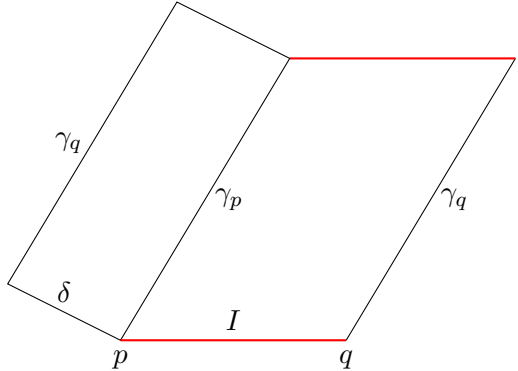
\begin{center}
\begin{tikzpicture}[scale=0.75]
\draw (4,0) -- (7,5); \draw[thick, red] (7,5) -- (3,5); \draw (3,5) -- (1,6) -- (-2,1) -- (0,0) -- cycle; \draw[thick, red] (0,0) -- (4,0);

\draw (0,0) node[below] {$p$}; \draw (4,0) node[below] {$q$}; \draw (5.5,2.5) node[right] {$\gamma_q$}; \draw (1.5,2.5) node[right]
{$\gamma_p$}; \draw (2,0) node[above] {$I$}; \draw (-1,0.5) node[above] {$\delta$}; \draw (-0.5,3.5) node[left] {$\gamma_q$};
\end{tikzpicture}

\end{center}

The closed geodesics $\gamma_p$ and $\gamma_q$ cut $\Sig$ into two cylinders, the one which contains $I$ will be denoted by $C_1$, the other one
by $C_2$. Let $\delta$ be a geodesic segment joining $p$ and $q$ which is contained in $C_2$.\\

\noindent The complement in $\Sig$ of the set $I\cup \gamma_p\cup \gamma_q \cup \delta$ is the union of two open parallelograms. By an embedding
of these two parallelograms into $\R^2$ which sends $\xi$ onto the constant vertical vector field $(0,1)$, we can associate the complex numbers
$Z,z,w$ to $I,\gamma_p,\delta$  respectively with a choice of orientation for each of these segments. Recall that $I$ is already oriented, hence
$Z$ is well defined, we can choose the orientation of $\gamma_p,$ and $\delta$ so that:

$$\mathbf{Area}(C_1)=\mathrm{Im}(Z\overline{z}) >0 \text{ and } \mathbf{Area}(C_2)=\mathrm{Im}(z\overline{w})>0.$$

\noindent We define two functions $\eta_1,\eta_2$ on $\C^3$ by the following formulae

 $$\eta_1(Z,z,w)=\mathrm{Im}(Z\overline{z}), \; \eta_2(Z,z,w)=\mathrm{Im}(z\overline{w}).$$

\noindent Set

$$\mathcal{D}=\{(Z,z,w)\in \C^3\; | \; \eta_1(Z,z,w)>0, \eta_2(Z,z,w)>0\}.$$

\noindent Remark that, given $(Z,z,w)$ in $\mathcal{D}$, one can construct a flat torus with a marked segment by first constructing two
parallelograms in $\R^2$ from the pairs of complex numbers $(Z,z)$ and $(z,w)$, and then gluing these two parallelograms as shown in the above
figure. We then get a map:

$$\rho :\mathcal{D} \lra \METTor,$$

\noindent which is surjective and locally homeomorphic. The pull-back of the volume form $\mu_{\mathrm{Tr}}$ on $\mathcal{D}$ is equal to the
Lebesgue measure of $\C^3$ up to a multiplicative constant. Clearly, the pull-back of the energy function $\mathcal{F}^\mathrm{et}$ on $\METTor$
is the following function on $\mathcal{D}$

$$\hat{\mathcal{F}}(Z,z,w)=\exp(-|Z|^2-(\eta_1(Z,z,w)+\eta_2(Z,z,w))).$$

\medskip

We say that a triple $(\Sig,I,\xi)$ is in {\em special position} if either $I$ is parallel to $\xi$, or the trajectory $(\psi_t)_{t>0}(p)$,
where $(\psi_t)$ is the flow generated by $\xi$, returns to $p$ without meeting any other point of $I$. Let $\METTor^\mathrm{sp}$ denote the set
of triples in special position in $\METTor$. We have

\begin{lemma}\label{FVTorL2}
The set $\METTor^\mathrm{sp}$ is of measure $0$ with respect to $\mu_{\mathrm{Tr}}$.
\end{lemma}

\dem The lemma follows from the fact that  $\METTor^\mathrm{sp}$ is the image under $\rho$ of the set

$$\{(Z,z,w)\in \mathcal{D}: \mathrm{Re} (Z)=0 \text{ or } \mathrm{Re}(z)=0 \},$$

\noindent which is obviously of measure zero with respect to the Lebesgue measure of $\C^3$. \carre

Now, let $(\Sig,I,\xi)$ be an element in the complement of  $\METTor^\mathrm{sp}$. Let $(Z,x,w)$ be the complex numbers associated to $I,
\gamma_p,$ and $\delta$ as above. Set

$$A=\mathrm{Re}(Z), a=\mathrm{Re}(z), b=\mathrm{Re}(w) \text{ and } B=\mathrm{Im}(Z), x=\mathrm{Im}(z), y=\mathrm{Im}(w).$$

\noindent Since $\xi$ is not parallel to $I$, we can take the direction $\theta$ in the proof of Lemma \ref{FVTorL1} to be the one determined by
$\xi$. Suppose that $\gamma_p$ arises from this construction then we have

\begin{equation*}\label{TorInEq1}
|a|\leq |A|.
\end{equation*}

\noindent Remark that, since $(\Sig,I,\xi)$ is not in special position, we have $|a|>0$. Since $C_2$ is a cylinder, we can choose the segment
$\delta$ so that

\begin{equation*}\label{TorInEq2}
|b|\leq |a|.
\end{equation*}

Now, set

$$\mathcal{D}_0=\{(Z,z,w) \in \mathcal{D} \; : \;  |A|\geq|a|\geq|b|\}.$$

\noindent From  the arguments above, we deduce that $\rho(\mathcal{D}_0)$ contains the complement of $\METTor^\mathrm{sp}$. Hence, the result of
Theorem \ref{FVTh} for this case will follow from the following proposition:

\begin{proposition} \label{FVTorPr} We have

$$\mathcal{J}=\int_{\mathcal{D}_0}\hat{\mathcal{F}}(Z,z,w)dAdBdadbdxdy=\int_{\mathcal{D}_0}\exp(-(A^2+B^2)-(\eta_1+\eta_2))dAdBdadbdxdy <\infty.$$
\end{proposition}

\dem From the definition of the domain $\mathcal{D}_0$, we have

$$\mathcal{J}= \int\int \exp(-(A^2+B^2))\times [\int_{-|A|}^{|A|}[\int_{-|a|}^{|a|} [\int\int \exp(-\eta_1-\eta_2)dxdy]db]da]dAdB.$$

\noindent Fix $A,B,a,b$, and consider the integral

$$\int\int \exp(-\eta_1-\eta_2)dxdy.$$

\noindent By definition, we have:

$$ \eta_1=Ba-Ax \text{ and } \eta_2= xb-ay.$$

\noindent Using the change of variables $(x,y)\longmapsto (\eta_1,\eta_2)$, we have

$$d\eta_1d\eta_2=|Aa|dxdy.$$

\noindent Since $\eta_1(Z,z,w)>0,$ and $\eta_2(Z,z,w)>0$ for every $(Z,z,w)$ in $\mathcal{D}_0$, it follows

$$\int\int_{(Z,z,w)\in
\mathcal{D}_0}\exp(-\eta_1-\eta_2)dxdy=\int_0^{+\infty}\int_0^{+\infty}\frac{e^{-\eta_1}e^{-\eta_2}}{|Aa|}d\eta_1d\eta_2 = \frac{1}{|Aa|}.$$

\noindent Consequently,

$$\mathcal{J}=
\int\int\exp(-A^2-B^2)[\int_{-|A|}^{|A|}[\int_{-|a|}^{|a|}\frac{1}{|Aa|}db]da]dAdB=4\int_{-\infty}^{\infty}\int_{-\infty}^{\infty}e^{-(A^2+B^2)}
dAdB < \infty.$$

\noindent This proves the proposition, and hence, Theorem \ref{FVTh} is proved for the case of $\METTor$.\carre

\section{Proof of Theorem \ref{FVTh}}\label{FVThPrf}

In this section, we will give the proof of Theorem \ref{FVTh} for the general case. Our strategy is very similar to the one in the particular
case $\METTor$, namely, we specify a finite family of open subsets of $\MET$ which covers a subset of full measure, and show that the integral
of the function $\mathcal{F}^{\mathrm{et}}$ on every member of this family is finite. Those open  subsets of $\MET$ are defined by means of
special admissible triangulations of surfaces in $\MET$ which are constructed by using the parallel vector field. Throughout this section, we
assume that $m<n$, which means that the family $\Ah=\{\A_1,\dots,\A_m\}$ contains at least a non-trivial tree. Note that the total number of
edges of the trees in $\Ah$ is $n-m$.\\

\subsection{Admissible matrix}

\noindent Set $N_2^*=N_2+(n-m)$, and $N=\dim_\C \MET$. Recall that we have

$$ N=\left\{%
\begin{array}{ll}
    N_1-N_2^*+1 & \hbox{ if $\alpha_i \in 2\pi\N, \; \forall i=1,\dots,n$}, \\
    N_1-N_2^* & \hbox{ otherwise}. \\
\end{array}%
\right.$$

\noindent Let us  define

\begin{definition}\label{AdmMtx}

A matrix $\mathbf{A}$ in $\mathbf{M}_{N_2^*,N_1}(\C)$ is called {\em admissible} if there exists an element $(\Sig,\hat{A},\xi)$ in $\MET$, and
an admissible triangulation $\T$ of $\Sig$ such that $\mathbf{A}$ is the coefficient matrix of the linear system associated to $\T$.\\

Let $a$ be a row of an admissible matrix. If $a$ corresponds to a triangle equation, then $a$ is called an {\em ordinary row}, otherwise, \ie
when $a$ corresponds to a boundary equation, it is called an {\em exceptional row}.\\

\end{definition}

Observe that the set of admissible matrices is finite. To see this, let $(\Sig,\hat{A},\xi)$ be an element of $\MET$, $\T$ be an admissible
triangulation of $\Sig$, and $\ST$ be the system associated to $\T$. Recall that $\ST$ consists of $N_2$ triangle equations, and $(n-m)$
boundary equations. Let $\AT \in \mathbf{M}_{N_2^*,N_1}(\C)$ be the coefficient matrix of $\ST$. Let $a$ be a row vector of $\AT$, then either

\begin{itemize}
\item[.] $a$ is an ordinary row, in this case, $a$ contains exactly three non-zero entries which belong to $\{\pm 1\}$, or

\item[.] $a$ is an exceptional row, in this case $a$ contains exactly two non-zero entries, one of which belongs to $\{\pm 1\}$, the other is of
the form $\pm e^{\imath\theta}$.

\end{itemize}

\noindent For any exceptional row, the angle $\theta$ belongs to a finite set of $[0;2\pi]$, since it corresponds to an edge of a tree the
forest $\hat{A}$, and is determined up to sign by the angles in $\underline{\alpha}$. As a consequence, we see that $a$ belongs to finite set of
$\C^{N_1}$. Therefore, $\AT$ belongs to a finite set of $\mathbf{M}_{N^*_2,N_1}(\C)$.\\

\medskip

Let $a$ be an exceptional row of an admissible matrix which is associated to an equation of the form

$$ \pm z_i \pm e^{\imath\theta} z_j=0.$$

\noindent We will call the operation consisting of multiplying $a$ by $e^{-\imath\theta}$ a {\em reversing operation}. Recall that $a$
corresponds to an edge of an erasing forest on a flat surface, and the angle $\theta$ is the rotation angle of the holonomy of a closed curve
which intersects the erasing forest at only one point in the corresponding edge transversely. Reversing the orientation of the closed curve
gives rise to the reversing operation on the row $a$.\\


Let $(\Sig,\hat{A},\xi)$ be an element of $\MET$. An admissible triangulation $\T$ of $\Sig$ does not give rise to a unique admissible matrix,
since the coefficients of the system $\ST$ depend on the following data

\begin{itemize}

\item[.] a numbering on the set of edges of the triangulation $\hat{\T}$, which is the triangulation induced by $\T$ on the surface obtained by
slitting open  $\Sig$ along trees in $\hat{A}$.

\item[.] a choice of orientation for each edge of $\hat{\T}$.

\item[.] a numbering on the set of triangles of $\hat{\T}$.

\item[.] a choice of orientation for the boundary of each triangle of $\hat{\T}$.

\item[.] for each edge of the forest $\hat{A}$, a choice of orientation for the closed curve which intersects $\hat{A}$ at only one point in
this edge transversely.

\end{itemize}

\noindent Therefore, we have an equivalence relation on the set of admissible matrices defined as follows

\begin{definition}

Two admissible matrices $\mathbf{A}_1$ and $\mathbf{A}_2$ are said to be equivalent if $\mathbf{A}_2$ can be obtained from $\mathbf{A}_1$ by a
sequence of the following operations

\begin{itemize}
\item[$\bullet$] interchanging two columns,

\item[$\bullet$] interchanging two rows,

\item[$\bullet$] changing sign of a columns,

\item[$\bullet$] changing sign of a row,

\item[$\bullet$] reversing operation on an exceptional row.
\end{itemize}

\end{definition}

\noindent Clearly, two admissible matrices arising from the same admissible triangulation are equivalent.\\

Let $\Ad$ denote the set of equivalence classes of admissible matrices in $\mathbf{M}_{N_2^*,N_1}(\C)$. For each $s$ in $\Ad$, choose a matrix
$\mathbf{A}_s$ in the equivalence class $s$, we then get a finite family $\{\mathbf{A}_s, \; s\in \Ad\}$ of matrices in
$\mathbf{M}_{N_2^*,N_1}(\C)$. We will associate to each $s$ in $\Ad$ an open subset of $\ker \mathbf{A}_s$ on which one can
define a map $\Phi_s$ with image in $\MET$ which is locally homeomorphic.\\

\noindent Given $s$ in $\Ad$, for any $Z=(z_1,\dots,z_{N_1})$ in $\ker \As$, such that $z_i\neq 0$, for $i=1,\dots,N_1$, let $\Sig_Z$ denote the
`surface' obtained from $Z$ by the following construction

\begin{itemize}

\item[1.] Construct a triangle in $\R^2$ from $z_i,z_j,z_k$ whenever there is an ordinary row $a$ in $\mathbf{A}_s$ such that

$$\DS{a\cdot \transpose{Z}= \pm z_i \pm z_j \pm z_k}.$$

\item[2.] Glue the triangles obtained from 1. together by identifying sides corresponding to the same coordinate of $Z$.

\item[3.] Identify the sides corresponding to $z_i$ and $z_j$ whenever there exists an exceptional row $a$  in $\mathbf{A}_s$ such that

$$\DS{a\cdot \transpose{Z}= \pm z_i \pm e^{\imath\theta}z_j}.$$

\end{itemize}

\noindent Let $\U_s$ be the open subset of $\ker \mathbf{A}_s$ which is defined by the condition:

\begin{tabbing}
  $\U_s=$ \= $\{ Z $ in $\ker \mathbf{A}_s$  with non-zero coordinates, such that $\Sig_Z$ is a closed, oriented, connected \\
          \>   flat surface, having exactly $n$ singularities with cone angles  $\alpha_1,\dots,\alpha_n \}$.
\end{tabbing}

\noindent We can then define a map $\Phi_s$ from $\U_s$ into $\MET$ by associating to a vector $Z$ in $\U_s$ the triple
$(\Sig_Z,\hat{A}_Z,\xi_Z)$, where $\hat{A}_Z$ is the forest consisting of the segments arising from the exceptional rows in $\mathbf{A}_s$, and
$\xi_Z$ is the vector field corresponding to the vertical constant vector field $(0,1)$ of $\R^2$.\\

\noindent By construction, for any point $(\Sig,\hat{A},\xi)$ in $\Phi_s(\U_s)$, there is an admissible triangulation $\T$ of $\Sig$ such that
the a local chart $\Psi_\T$ defined in a neighborhood of $(\Sig,\hat{A},\xi)$ verifies $\Phi^{-1}_s=\Psi_\T$. It follows that $\Phi_s(\U_s)$ is
an open subset of $\MET$. Since every element of $\MET$ is contained in the domain of a local chart associated to an admissible triangulation,
the following proposition is now clear

\begin{proposition}\label{FVolETpr1}
The family $\{\Phi_s(\U_s), s\in \Ad\}$ is a finite open cover of  the space $\MET$.
\end{proposition}

\subsection{Primary, auxiliary systems of indices, and admissible triple}

Set

$$K=N-(2g+m-2)=\left\{%
\begin{array}{ll}
    n-m+1 & \hbox{ if $N=2g+n-1$,} \\
    n-m & \hbox{if $N=2g+n-2$.} \\
\end{array}%
\right.$$

\noindent In what follows, we will identify any matrix in $\mathbf{M}_{N_2^*,N_1}(\C)$ (resp. $\mathbf{M}_{N_2,N_1}(\C))$ to the linear map from
$\C^{N_1}$  to $\C^{N^*_2}$ (resp. to $\C^{N_2}$) which is defined by this matrix in the canonical bases of $\C^{N_1}$, and $\C^{N^*_2}$ (resp.
of $\C^{N_1}$, and $\C^{N_2}$).\\

\begin{definition}\label{DefPrSyst}
Given a matrix $\mathbf{A}$ in $\mathbf{M}_{n_1,n_2}(\C)$ with $n_1<n_2$, set $r=\dim \ker\mathbf{A}= n_2-\mathrm{rk} \mathbf{A}$. A {\em
primary system of indices} for $\mathbf{A}$ is an ordered subset $(i_1,\dots,i_{r})$ of $\{1,\dots,n_2\}$ such that there exist $n_2$ complex
linear functions ${f_i :\C^{r}\ra \C, \; i=1,\dots,n_2}$, verifying the following condition:

$$ (z_1,\dots,z_{n_2}) \in \ker \mathbf{A} \text{ if and only if } z_i=f_i(z_{i_1},\dots,z_{i_{r}}), \text{ for } i=1,\dots,n_2.$$
\end{definition}


\begin{definition}\label{DefAuSyst}

Given an $s$ in $\Ad$, and a primary system of indices $I=(i_1,\dots,i_{N})$ for $\mathbf{A}_s$, an {\em auxiliary system of indices} for $I$ is
an ordered subset $(j_K,\dots,j_{N})$ of $\{1,\dots,N_1\}$, which is empty if $K>N$ (that is when $g=0$, and $m=1$), such that, for
$k=K,\dots,N$

\begin{itemize}

\item[$i)$] $f_{j_k}$ depends only on $(z_{i_1},\dots,z_{i_{k-1}})$,

\item[$ii)$] the coefficients of $z_{i_K},\dots,z_{i_{k-1}}$ in $f_{j_k}$ are all real,

\item[$iii)$] There exists an ordinary row in $\mathbf{A}_s$ whose $i_k$-th and $j_k$-th entries are both non-zero.

\end{itemize}

\end{definition}

\noindent \underline{\bf Convention:} Given a matrix $\mathbf{A}$ in $\mathbf{M}_{N_2^*,N_1}(\C)$, or in $\mathbf{M}_{N_2,N_1}(\C)$, in what
follows, we will say that $z_j$ is a linear function of $(z_{i_1},\dots,z_{i_k})$, or $z_j$ depends linearly on $(z_{i_1},\dots,z_{i_k})$ as
$(z_1,\dots,z_{N_1})$ varies in $\ker \mathbf{A}$ if there exists a vector $(\lbd_1,\dots,\lbd_k)$ in $\C^{k}$ such that

$$\mathbf{A}\cdot ^t(z_1,\dots,z_{N_1})=0 \text{ implies } z_j=\lbd_1z_{i_1}+\dots+\lbd_{k}z_{i_k}.$$

\rem  If $(j_K,\dots,j_N)$ is an auxiliary system for $(i_1,\dots,i_N)$, then  we have

\begin{itemize}
\item[$\bullet$] $z_{j_k}$ can be written as a linear function of $(z_{i_1},\dots,z_{i_{k-1}})$, for  $k=K,\dots,N$, as $Z=(z_1,\dots,z_{N_1})$
varies in $\ker\As$.

\item[$\bullet$] Assume that $(\Sig,\hat{A},\xi)=\Phi_s(Z)$, and let $\T$ be the geodesic triangulation of $\Sig$ which is obtained from the
construction of $\Phi_s$, then the condition $iii)$ of \ref{DefAuSyst} implies that $z_{i_k}$ and $z_{j_k}$ are associated to two sides of a
triangle in $\T$.

\end{itemize}

For the case $\METTor$, let  $(\Sig, I,\xi)$ be an element of $\METTor$, and let $p,q,\gamma_p,\gamma_q,\delta$, and $Z,z,w$ be as in Section
\ref{TorSect}. We can add some geodesic segments whose endpoints are contained in the set $\{p,q\}$ to get a triangulation of $\Sig$. We then
get a triangulation of the surface $\hat{\Sig}$ which is obtained from $\Sig$ by slitting along $I$. This triangulation gives rise to a an
admissible matrix $\mathbf{A}$ in $\mathbf{M}_{5,7}(\C)$ with $\dim \ker\mathbf{A}=3$. There exists a linear isomorphism $\varphi$ from $\C^3$
to $\ker\mathbf{A}$, we can arrange so that, for any $(z_1,\dots,z_7)=\varphi(Z,z,w)$ then $z_1=Z, z_2=z, z_3=w$. In this case, $N=3, K=2$,
therefore $(1,2,3)$ is a primary system of indices for $\mathbf{A}$, and $(1,2)$ is an auxiliary system of indices for $(1,2,3)$.\\

By Proposition \ref{FVolETpr1}, we know that $\MET$ is covered by the family of open subsets $\{\Phi_s(\mathcal{U}_s), \; s \in \Ad\}$.
Therefore, to prove Theorem \ref{FVTh}, we only need to show that the integral of the function $\mathcal{F}^\mathrm{et}$ on $\Phi_s(\U_s)$ is
finite. This would be done if we could show that the integral of $\Phi_s^*\mathcal{F}^\mathrm{et}$ on $\U_s$ is finite. However, the domain
$\U_s$ is still too large, and this integral can be infinite. The primary and auxiliary systems of indices for $\As, \; s \in \Ad,$ will allow
us to specify a finite family of sub-domains of $\U_s$ on which the integral of $\Phi_s^*\mathcal{F}^\mathrm{et}$ is finite, and whose images
under $\Phi_s$ cover a full measure subset of $\Phi_s(\U_s)$.\\

Consider $\As$ for some $s$ in $\Ad$. Let $a_1,\dots,a_{N_2}$ denote the ordinary rows, and $b_1,\dots,b_{n-m}$ denote the exceptional rows of
$\mathbf{A}_s$. Let $\mathbf{A}^*_s \in \mathbf{M}_{N_2,N_1}(\C)$ be the matrix consisting of the ordinary rows of $\mathbf{A}_s$, and set
$\tilde{N}=\dim \ker \mathbf{A}^*_s$.\\

\begin{definition}\label{DefBdrInd}
If the $i$-th column of $\mathbf{A}^*_s$ has only one non-zero entry, we say that $i$ is a {\em boundary index} of $\As$. Two boundary indices
$i_1$ and $i_2$ are said to be {\em paired up} if there exists an exceptional row in $\mathbf{A}_s$ whose $i_1$-th and $i_2$-th entries are
non-zero.
\end{definition}

Fix a vector $Z=(z_1,\dots,z_{N_1})$ in $\U_s$, and  let $(\Sig,\hat{A},\xi)$ be the image of $Z$ under $\Phi_s$. Recall that $\Sig$ comes along
with an admissible triangulation $\T$. Let $\hSig$ denote the surface obtained by slitting open $\Sig$ along $\hat{A}$, and $\hT$ denote the
triangulation of $\hSig$ which is induced by $\T$. By definition, the coordinates of $Z$ is in bijection with the set of edges of a
triangulation $\hT$, and the rows of $\Asv$ is in bijection with the set of triangles of $\hT$ . If $i$ is a boundary index of $\As$, then $z_i$
corresponds to an edge of $\hT$ which is contained in the boundary of $\hSig$. If $i_1$, $i_2$ are two boundary indices which are paired up,
then the edges corresponding to $z_{i_1}$, and $z_{i_2}$ arise from the same edge of a tree in the forest $\hat{A}$. Observe that the set of
boundary indices of $\As$ contains exactly $2(n-m)$ elements divided into $(n-m)$ pairs. First, let us prove the following

\begin{lemma}\label{FVlm1bis}

We have $\rk\Asv=N_2$, or equivalently $\tilde{N}=N_1-N_2=(2g+n-2)+(n-m)$.

\end{lemma}

\dem We will show that the row vectors $(a_1,\dots,a_{N_2})$ are linearly independent in $\C^{N_1}$. Assume that there exists
$(\lbd_1,\dots,\lbd_{N_2})$ in $\C^{N_2}$ such that

$$\lbd_1 a_1 +\dots+\lbd_{N_2}a_{N_2}=0.$$

\noindent Observe that, if $a_{i_1}$ and $a_{i_2}$ correspond to two adjacent triangles of $\hT$, then there exists $j\in \{1,\dots,N_1\}$ such
that the $j$-th column of $\Asv$ contains exactly two non-zero entries, on the $i_1$-th and  the $i_2$-th rows. It follows that if
$\lbd_{i_1}=0$, then $\lbd_{i_2}=0$.\\

\noindent Since $n-m>0$, the set of boundary indices is non-empty, which means that there exists a column in $\Asv$ which contains exactly one
non-zero entry. Hence, there exists $j \in \{1,\dots,N_2\}$ such that $\lbd_j=0$. Since the surface $\hSig$ is connected, the argument above
implies that $\lbd_1=\dots=\lbd_{N_2}=0$, and the lemma follows. \carre

The next lemma tells us that a primary system of indices for $\Asv$ contains at most $2(n-m)-1$ boundary indices.

\begin{lemma}\label{FVlm2bis}
All the boundary indices can not be contained in a primary system of indices for $\Asv$.
\end{lemma}

\dem Recall that the sign of a row in $\Asv$ is determined by a choice of orientation on the boundary of the corresponding triangle in $\hT$.
Note that we are free to permute, and change sign of rows and columns in $\Asv$. Let $I_{b} \subset \{1,\dots,N_1\}$ be the subset of boundary
indices for $\Asv$.\\

\noindent For each triangle of $\hT$, we choose the orientation of its boundary coherently with the orientation of the surface. Since each edge
in the interior of $\hSig$ belongs to two distinct triangles, it follows that we have

$$(a_1+\dots+a_{N_2})\cdot \transpose{Z}=\sum_{i\in I_b} \pm z_i.$$

\noindent Therefore, for every $(z_1,\dots,z_{N_1})$ in $\ker \Asv$, we have

\begin{equation}\label{Eq0}
 \sum_{i\in I_b}\pm z_i=0
\end{equation}

\noindent which implies that the family of coordinates $(z_i, \; i \in I_b)$ is not linearly independent as $(z_1,\dots,z_{N_1})$ varies in
$\ker\Asv$, and the lemma follows. \carre

Our goal now is to prove that there exists a primary system of indices for $ \mathbf{A}_s$ whose $(K-1)$ first indices are boundary indices. Let
us first prove the following

\begin{lemma}\label{FVlm1}
There exist primary systems of indices for $\Asv$ whose first $2(n-m)-1$ elements are boundary indices.
\end{lemma}

\dem Assume that $I_b=\{1,\dots,2(n-m)\}$ is the set of boundary indices of $\Asv$. By permuting the rows of $\Asv$, we can assume that the only
non-zero of the first column is on the first row, that is the first entry of $a_1$ is $\pm 1$. We will show that, as $(z_1,\dots,z_{N_1})$
varies in $\ker \Asv$, the coordinates $(z_2,\dots,z_{2(n-m)})$ are linearly independent, that is $\ker \Asv$ is not contained in the kernel  of
any any linear function of the  form

$$f: (z_1,\dots,z_{N_1}) \mapsto \lambda_2z_2+\dots+\lambda_{2(n-m)}z_{2(n-m)}.$$

\noindent It follows that we can add $(\tilde{N}-2(m-n)+1)$ indices to the family $\{2,\dots,2(n-m)\}$  to get a primary system of indices for
$\Asv$, which proves the lemma.\\

All we have to show is that, if there exists a vector $\lambda'=(\lambda'_1,\dots,\lambda'_{N_2})\in \C^{N_2}$ such that

\begin{equation}\label{Eq1}
(\sum_{i=1}^{N_2}\lambda'_i a_i)\cdot \transpose{(z_1,\dots,z_{N_1})}=\sum_{i=2}^{2(n-m)}\lambda_iz_i
\end{equation}

\noindent then $\lambda'_i=0, \; i=1,\dots,N_2$.\\

\noindent First, observe that we must have $\lambda'_1=0$, since there is only one non-zero entry in the first column of $\Asv$. Consider two
adjacent triangles $\Delta_1$, $\Delta_2$ of $\hT$. Each common edge of $\Delta_1$ and $\Delta_2$ correspond to a coordinate $z_j$ of $Z$, with
$j>2(n-m)$. Let $a_{i_1},a_{i_2}$ be the rows in $\Asv$ which correspond to $\Delta_1, \Delta_2$ respectively, then, in the $j$-th column of
$\Asv$ there exactly two non-zero entries, on the $i_1$-th and the $i_2$-th rows. Now, since the right hand side of (\ref{Eq1}) does not contain
$z_j$, with $j>2(n-m)$, we deduce that, if $\lambda'_{i_1}=0$, then $\lambda'_{i_2}=0$.  We already have $\lambda'_1=0$, and  since $\hat{\Sig}$
is connected, it follows that $\lambda'_i=0$, for  $i=1,\dots,N_2$. \carre

\begin{lemma}\label{FVlm1b}
There exist primary systems of indices for $\As$ whose first $(K-1)$ elements are boundary indices.

\end{lemma}

\dem We can assume that the set of boundary indices of $\Asv$ is $\{1,\dots,2(n-m)\}$, and that $i$ and $(n-m)+i$, $i=1,\dots,n-m$, are paired
up, which means that any $(z_1,\dots,z_{N_1})$ in $\ker\As$ satisfies $(n-m)$ equations of the form

\begin{equation}\label{Eq6}
z_i\pm e^{\imath\theta_i}z_{(n-m)+i}=0, \; i=1,\dots,n-m,
\end{equation}

\noindent with some $\theta_i$ in a finite set. Since $(z_1,\dots,z_{N_1})$ also satisfies the equation (\ref{Eq0}), it follows that we have

\begin{equation}\label{Eq7}
\sum_{i=1}^{n-m}(1\pm e^{\imath\theta_i})z_i=0
\end{equation}

By Lemma \ref{FVlm1}, we know that there exists a primary system of indices $\tilde{I}$ for $\Asv$ whose $2(n-m)-1$ first elements are boundary
indices. We will show that a primary system of indices for $\As$ can be obtained by removing some boundary indices in $\tilde{I}$. We have two
issues:

\begin{itemize}
\item[$\bullet$] \underline{Case 1:} $\alpha_i\in 2\pi\N, \; i=1,\dots,n$. In this case, $N=2g+n-1$, $K=(n-m)+1$, and $\tilde{N}=(2g+n-2)+(n-m)=
N+(n-m)-1$. Note that in this case, all the angles $\theta_i$ are zero, and with  appropriate choices of orientation for the edges of $\hT$ in
the boundary of $\hSig$, the equation (\ref{Eq7}) is trivial (cf. \cite{Ng1}).\\

Let $I$ be the ordered subset of $\{1,\dots,N_1\}$ which is obtained by removing the indices $\{(n-m)+1,\dots,2(n-m)-1\}$ from $\tilde{I}$. The
set $I$ contains $n-m=K-1$ boundary indices. Let us show that $I$ is a primary system of indices for $\As$. First, observe that, for any
$(z_1,\dots,z_{N_1})$ in $\ker\As$, $z_i, \; i=1,\dots,N_1$ can be written as a linear function of $\{z_k, \; k\in \tilde{I}\}$, since
$\tilde{I}$ is a primary system of indices for $\Asv$. Using the equations (\ref{Eq6}), we can replace $z_{(n-m)+j}$ by $\pm
e^{\imath\theta_j}z_j, \; j=1,\dots,n-m$. Therefore, $z_i, \; i=1,\dots,N_1$, can be written as a linear function of $(z_k, \; k\in I)$ as
$(z_1,\dots,z_{N_1})$ varies in $\ker\As$. Moreover, we have

$$\card\{I\}=2g+n-1=\dim \ker\As,$$

which implies that $I$ is a primary system of indices for $\As$.\\

\item[$\bullet$] \underline{Case 2:} there exist $i\in \{1,\dots,n\}$ such that $\alpha_i\not\in 2\pi\N$. In this case, $N=2g+n-2$, $K=n-m$,
$\tilde{N}-N=n-m$, and the equation (\ref{Eq7}) is non-trivial (cf. \cite{Ng1}). Without loss of generality, we can assume that the coefficient
of $z_1$ in (\ref{Eq7}) is non-zero, which means that, as $(z_1,\dots,z_{N_1})$ varies in $\ker \As$,  $z_1$ can be  written as a linear
function of $(z_2,\dots,z_{n-m})$.\\

Let $I$ be the ordered subset of $\{1,\dots,N_1\}$ which is obtained by removing the indices $\{1, (n-m)+1,\dots,2(n-m)-1\}$ from $\tilde{I}$.
Clearly, the set $I$ contains $(n-m)-1=K-1$ boundary indices. Since $\tilde{I}$ is a primary system of indices for $\Asv$, using the equations
(\ref{Eq6}), and (\ref{Eq7}), we see that, as $(z_1,\dots,z_{N_1})$ varies in $\ker\As$,  for $i=1,\dots,N_1$, $z_i$ can be written as a linear
function of $(z_k, \; k\in I)$. Moreover, we have

$$\card\{I\}=N=\dim \ker\As,$$

therefore, $I$ is a primary system of indices for $\As$. The proof of the lemma is now complete. \carre

\end{itemize}

We can now define

\begin{definition}\label{DefAdmTrp}
Any triple $(\As;I;J)$, where $I$ is a primary system of indices for $\As$ whose $(K-1)$ first elements are boundary indices, and $J$
is an auxiliary system od indices for $I$, will be called an {\em admissible triple}.\\
\end{definition}

\noindent Clearly, the number of admissible triples is finite.\\


\subsection{Good triangulation} \label{GTrSbSect}
Throughout this subsection, given a point $(\Sig,\hat{A},\xi)$ in $\MET$, we denote by $\hSig$ the flat surface with piecewise geodesic boundary
obtained by slitting open $\Sig$ along the trees in $\hat{A}$. Let $\hat{V}$ denote the finite subset of $\hSig$ which arises from the vertex
set $V$ of $\hat{A}$. The vector field $\xi$ of $\Sig$ gives rise to a parallel vector field of $\hSig$ which will be denoted again by $\xi$.
For any admissible triangulation $\T$ of $(\Sig,\hat{A},\xi)$, let $\hT$ denote the induced triangulation of $\hSig$.\\

Let $(\psi_t), \; t\in \R,$ denote the flow generated by $\xi$ on $\hSig$. Given  a point $p$ in $\inter(\hSig)\setm \hat{V}$, if there exists
$t_0>0$ (resp. $t_0<0$) such that $\psi_{t_0}(p)\in \hat{V}\cup\partial \hSig$, then, for every $t>t_0$ (resp. $t<t_0$), we consider, by
convention, that $\psi_t(p)=\psi_{t_0}(p)$. \\

Let $a$ be a geodesic segment contained in the boundary of $\hSig$ with endpoints in $\hat{V}$. We can extend the field $\xi$ by continuity to
$\inter(a)$. Assume that $a$ is not parallel to the field $\xi$, then we say that $a$ is an {\em upper} (resp. {\em lower}) boundary segment, if
the field $\xi$ on $\inter(a)$ points outward (resp. inward). Observe that in this case, the image of $\inter(a)$ by $\psi_t$ is well defined
for all $t\in \R$.\\

Let $(\Sig,\hat{A},\xi)$ be a point in $\MET$, and let $e$ be a geodesic segment of $\hSig$ with endpoints in $\hat{V}$, we denote by $h(e)$ the
{\em transversal measure} of $e$ with respect to $\xi$ which is defined as follows: if we choose an isometric embedding of a neighborhood of $e$
into $\R^2$ such that the vector field $\xi$ is mapped to the constant vertical vector field $(0,1)$ of $\R^2$, then $h(e)$ is  the length of
the horizontal projection of the image of $e$. We call $h(e)$ the {\em horizontal length} of $e$.\\

\noindent A triangle in $\hSig$ whose sides are geodesic segments denoted by $e_1,e_2,e_3$ is said to be {\em good} if $h(e_i)>0$, for
$i=1,2,3$. Given a good triangle $\Delta$, we call the unique side of $\Delta$ of maximal horizontal length the {\em base} of $\Delta$. Let
$\hT$ be a triangulation of $\hSig$ which arises from an admissible triangulation of $\Sig$,  if all of triangles of $\hT$ are good, then $\hT$
is called a {\em good triangulation}. The following proposition asserts that a `generic' element always admits a good triangulation.\\

\begin{proposition} \label{FVpr2}

Let $(\Sig, \hat{A},\xi)$ be an element of $\MET$. Suppose that there exist no geodesic segments in $\hSig$ with endpoints in $\hat{V}$ which
are parallel to the field $\xi$, then there exists a good triangulation $\hT$ of $\hSig$ whose edges are denoted by $\{e_1,\dots,e_{N_1}\}$ so
that,

\begin{itemize}

\item[$\bullet$] The edges of $\hT$ in the boundary of $\hSig$ are denoted by $\{e_1,\dots,e_{2(n-m)}\}$.

\item[$\bullet$] For every $i\in\{2(n-m)+1,\dots,N_1\}$, there exists $j<i$, and a triangle $\Delta$ of $\hT$ whose boundary contains both
$e_i,e_j$ such that $e_j$ is the base of $\Delta$.

\end{itemize}

\end{proposition}

\dem We construct a geodesic triangulation of $\hSig$ whose vertex set is $\hat{V}$ as follows: let $e_1,\dots,e_{2(n-m)}$ denote the geodesic
segments in the boundary of $\hSig$ with endpoints in $\hat{V}$. Assume that the segment $e_{1}$ is of maximal horizontal length among the set
$\{e_1,\dots,e_{2(n-m)}\}$. By assumption, we have $h(e_{1})>0$. Let $p,q$ denote the endpoints of $e_{1}$ (it may happen that $p\equiv q$).
Consider the following procedure:\\

\noindent Assume that $e_{1}$ is a lower boundary segment, consider the stripe $S_t$ swept by $\{\psi_t(\inter(e_1)), t>0\}$. Since
$h(e_{1})>0$, for some $t$ finite, this stripe must meet a point in the set $\hat{V}\cup \partial\hSig$, otherwise its area would tend to
infinity as $t$ tends to $+\infty$.\\

\noindent Since the horizontal length of $e_{1}$ is maximal among the set $\{h(e_1),\dots,h(e_{2(n-m)})\}$, suppose that, for some $t\in \R^+$,
$\psi_t(\inter(e_{1}))$ is contained in a the geodesic segments $e_i$ in the set $\{e_1,\dots, e_{2(n-m)}\}$, then we must have
$\overline{\psi_t(\inter(e_1))}=e_i$. This implies that there is a geodesic segment parallel to the field $\xi$ joining $p$ to a point in
$\hat{V}$, which is a contradiction to the assumption of the lemma. Therefore, there exists a smallest value $t_0>0$ such that
$\psi_{t_0}(\inter(e_{1}))$ contains a point in $\hat{V}$.\\

\noindent  Let $r$ be a point in $\psi_{t_0}(\inter(e_{1}))\cap \hat{V}$, and let $e', e''$ denote the two geodesic segments contained in the
stripe $S_{t_0}$ which join $r$ to $p$, and to $q$.  Note that even though $p$ and $q$ may coincide, the two segments $e'$, and $e''$ are always
distinct. It can happen that one of the segments $e',e''$ already appears in the set $\{e_1,\dots,e_{2(n-m)}\}$ but not both of them, unless
$\hSig$ is a triangle. By assumption, we have $h(e')>0$, and $h(e'')>0$, and by construction, $e_{1}$ is the base of the good triangle bounded
by $e',e''$, and $e_{1}$. We will call $e_{1}$ the {\em supporter} of $e'$ and $e''$.\\


In the case where $e_{1}$ is an upper boundary segment, by considering $\{\psi_t(\inter(e_{1})), \; t<0\}$ instead of $\{\psi_t(\inter(e_{1})),
\; t>0\}$, we also get a similar result.\\

\noindent Cut off the triangle bounded by $e_{1}, e', e''$ from the surface $\hSig$ along the segments $e'$ and $e''$. The remaining surface is
a flat surface with piecewise geodesic boundary which is not necessarily connected. On this new surface, we still have a parallel vector field
which is the restriction of $\xi$. We can now reapply the same procedure to the new surface. The assumption of the proposition allows us to
continue this procedure until the surface $\hSig$ is cut into triangles with vertices in $\hat{V}$, that is until we get a geodesic
triangulation $\hT$ of $\hSig$ whose vertex set is $\hat{V}$, this triangulation is necessarily  a good triangulation.\\

We number the edges of $\hT$ which are contained in the interior of $\hSig$ according to their appearing order in the  procedure above, the
ordering of two edges which appear in the same step is not important. Since every edge of $\hT$ in the interior of $\hSig$ admits a supporter
which appears in the procedure before itself, the proposition is then proved. \carre

\begin{proposition}\label{FVpr2b}

If $(\Sig,\hat{A},\xi)$ is a point in  $\MET$ satisfying the condition of Proposition \ref{FVpr2}, then there exists an admissible triple
$(\As;I;J)$, where  $I=(i_1,\dots,i_N)$, and $J=(j_K,\dots,j_N)$, and a vector $Z^0=(z^0_1,\dots,z^0_{N_1})$ in $\U_s$ such that

\begin{itemize}
\item[$\bullet$] $(\Sig,\hat{A},\xi)=\Phi_s(Z^0)$.

\item[$\bullet$] $|\mathrm{Re}(z^0_{j_k})| > |\mathrm{Re}(z^0_{i_k})|$ for any $k=K,\dots,N$.


\end{itemize}

\end{proposition}

\dem Let $\hT$ be the good triangulation of $\hSig$ which is obtained from Proposition \ref{FVpr2}. Let $\AT$ be the matrix in
$\mathbf{M}_{N_2^*,N_1}(\C)$ associated to $\hT$.  Let $Z^0=(z^0_1,\dots,z^0_{N_1})$ denote the vector of $\ker\AT$ whose coordinates are
associated to the edges of $\hT$. In what follows, we consider any vector $Z=(z_1,\dots,z_{N_1})$ in $\C^{N_1}$ as a function from the set of
edges of $\hT$ to $\C$ such that $z_i= Z(e_i)$.\\

\noindent By construction, the set $I_b$ of boundary indices for $\AT$ is $\{1,\dots,2(n-m)\}$. Let $\ATv$ be the matrix in
$\mathbf{M}_{N_2,N_1}(\C)$ consisting of all ordinary rows of $\AT$. Let $N$, and $\tilde{N}$ denote the dimensions of $\ker\AT$, and $\ker
\ATv$ respectively. We first choose a primary system of indices $\tilde{I}$ for $\ATv$ as follows:

\begin{itemize}
\item[$\bullet$] The first $2(n-m)-1$ elements of $\tilde{I}$ are $\{2,\dots,2(n-m)\}$, by Lemma \ref{FVlm1}, we know that, as
$Z=(z_1,\dots,z_{N_1})$ varies in $\ker \ATv$, the family of coordinates $(z_2,\dots,z_{2(n-m)})$ is linearly independent.

\item[$\bullet$] Assume that we have chosen $k$ indices $(i'_1,\dots,i'_k)$ for $\tilde{I}$, then the index $i'_{k+1}$ of $\tilde{I}$ is the
smallest index $i'$ such that, as $(z_1,\dots,z_{N_1})$ varies in $\ker\ATv$, $z_{i'}$ can not be written as a linear function of
$(z_{i'_1},\dots, z_{i'_k})$, in other words, the family of coordinates $(z_{i'_1},\dots,z_{i'_k},z_{i'})$ is linearly independent.

\end{itemize}

By Lemma \ref{FVlm2bis}, we know that, for $k=2(n-m),\dots,\tilde{N}$, $i'_k$ is not a boundary index, that is $i'_k > 2(n-m)$. For any $k$ in
$\{2(n-m),\dots,\tilde{N}\}$, consider the edge $e_{i'_k}$ of $\hT$. From Proposition \ref{FVpr2}, we know that there exists an edge $e_{j'_k}$
with $j'_k < i'_k$, and a triangle $\Delta_k$ of $\hT$ whose boundary contains both $e_{i'_k}$, and $e_{j'_k}$ such that $e_{j'_k}$ is the base
of $\Delta_k$. Consequently, we have

\begin{equation}\label{InEq0}
|\mathrm{Re}(z^0_{j'_k})|=h(e_{j'_k}) > h(e_{i'_k})=|\mathrm{Re}(z^0_{i'_k})|
\end{equation}

Let $J$ denote the ordered subset $(j'_{2(n-m)},\dots,j'_{\tilde{N}})$ of $\{1,\dots,N_1\}$. From the definition of $i'_k$, for
$k=2(n-m),\dots,\tilde{N}$, as $(z_1,\dots,z_{N_1})$ varies in $\ker\ATv$, we can write

$$z_{j'_k}=\tilde{f}_{j'_k}(z_{i'_1},\dots,z_{i'_{k-1}}),$$

\noindent where $\tilde{f}_{j'_k}$ is some fixed linear function. Since the matrix $\ATv$ is real, all the coefficients of $\tilde{f}_{j'_k}$
are also real.\\

By Lemma \ref{FVlm1b}, we know that, by removing $K'=2(n-m)-K$ boundary indices from $\tilde{I}$, we obtain a primary system of indices $I$ for
$\AT$ whose first $(K-1)$ elements are boundary indices. We will show that $J$ is an auxiliary of $I$, which, together with (\ref{InEq0}), will
allow us to conclude.\\

First, observe that we can write $I=(i_1,\dots,i_N)$, where $i_1,\dots,i_{K-1}$ are boundary indices, and for $k=K,\dots,N, \;  i_k=i'_{k+K'}$.
Since $\tilde{N}-N=2(n-m)-K=K'$, we can write $J=(j_{K},\dots,j_N)$, where $j_k=j'_{k+K'}$. As a consequence, for $k=K,\dots,N$, the condition
that there is a triangle in $\hT$ whose boundary contains both $e_{i_k}$, and $e_{j_k}$ is satisfied.\\

\noindent We already know that, as $(z_1,\dots,z_{N_1})$ varies in $\ker\AT \subset \ker \ATv$, for $k=K,\dots,N$, we have

$$z_{j_k}=z_{j'_{k+K'}}=\tilde{f}_{j_k}(z_{i'_1},\dots, z_{i'_{k+K'-1}}),$$

\noindent where $\tilde{f}_{j_k}$ is a linear function  with real coefficients. We can then transform $\tilde{f}_{j_k}$ into a linear function
$f_{j_k}$ of $(z_{i_1},\dots,z_{i_{k-1}})$ by using equations of the form (\ref{Eq6}), and (\ref{Eq7}). Since the equations (\ref{Eq6}), and
(\ref{Eq7}) involve only boundary indices, we deduce that the coefficients of $z_{i_K},\dots,z_{i_{k-1}}$ in $f_{j_k}$ are all real, which
allows us to conclude that $J$ is an auxiliary system of indices for $I$.\\

Clearly, the inequalities (\ref{InEq0}) can be rewritten as

\begin{equation}\label{InEq0b}
|\mathrm{Re}(z^0_{j_k})| >|\mathrm{Re}(z^0_{i_k})|, \; k=K,\dots,N
\end{equation}

\noindent We know that $\mathbf{A}_\T$ is equivalent to a matrix $\mathbf{A}_s$ with $s$ in $\Ad$. The transformation of $\mathbf{A}_\T$ into
$\mathbf{A}_s$ consists of renumbering the coordinates in $\C^{N_1}$, and changing their sign. By this transformation, $(i_1,\dots,i_N)$ becomes
a primary system of indices for $\As$, and $(j_K,\dots,j_N)$ becomes an auxiliary system of indices for $(i_1,\dots,i_N)$. Therefore, we get an
admissible triple $(\As;I;J)$, and a vector $Z^0$ in $\U_s$ which verify the conditions in the statement of the proposition. \carre

Now, given an admissible triple $(\As;I;J)$, where $I=(i_1,\dots,i_{N}), J=(j_K,\dots,j_N)$, let $\U_s(I;J)$ denote the following subset of
$\U_s$

$$\U_s(I;J)=\{(z_1,\dots,z_{N_1})\in \U_s\; : \; |\mathrm{Re}(z_{i_k})|\leq |\mathrm{Re}(z_{j_k})|,\;  k = K,\dots,N \}.$$

We say that the element $(\Sig,\hat{A},\xi)$ of $\MET$ is in {\em special position} if there exists a geodesic segment in $\hSig$ with endpoints
in $\hat{V}$ parallel to the field $\xi$. Let $\MET^\mathrm{sp}$ denote the subset of $\MET$ consisting of elements in special position. A
direct consequence of Proposition \ref{FVpr2b} is the following

\begin{corollary}\label{FVcor2}
The finite family $\DS{\{\Phi_s(\U_s(I;J)):\; (\mathbf{A}_s;I;J) \text{ is admissible }\}}$ covers the complement of $\MET^\mathrm{sp}$ in
$\MET$.
\end{corollary}

The next proposition tells us that $\MET\setm\MET^\mathrm{sp}$ is a subset of full measure in $\MET$.

\begin{proposition}\label{FVpr1}

 The set $\MET^\mathrm{sp}$ is a null set in $\MET$ with respect to $\mu_{\mathrm{Tr}}$.

\end{proposition}

\dem For every $s$ in $\Ad$, let $\mu_s$ denote the volume form on $\U_s$ which is the pull-back of $\mu_\Tr$ under $\Phi_s$. Let
$(\Sig,\hat{A},\xi)$ be a point in $\MET^\mathrm{sp}$, let $e$ be a geodesic segment of $\hSig$ with endpoint in $\hat{V}$ parallel to the field
$\xi$. There exists an admissible triangulation $\hT$ of $\hSig$ such that the $1$-skeleton of $\hT$ contains $e$. Since $e$ is parallel to
$\xi$, the complex number associated to $e$ in the local chart arising from $\hT$ is purely imaginary. As a consequence, there exist

\begin{itemize}

\item[.] $s\in\Ad$,

\item[.] $i\in \{1,\dots,N_1\}$, and

\item[.] $Z\in \{(z_1,\dots,z_{N_1})\in \U_s|\;  \mathrm{Re}(z_i)=0\}$,

\end{itemize}

\noindent such that $(\Sig,\hat{A},\xi)=\Phi_s(Z)$. For every $s\in\Ad$, and every $i\in \{1,\dots,N_1\}$, set

$$\U_s^i=\U_s\cap\{(z_1,\dots,z_{N_1})\in \C^{N_1} \;| \; \mathrm{Re}(z_i)=0\}.$$

\noindent Note that if $Z\in \U^i_s$, then $\Phi_s(Z)\in \MET^{\mathrm{sp}}$. It follows that

$$\MET^\mathrm{sp} =\bigcup_{s\in \Ad}\bigcup_{i=1}^{N_1}\Phi_s(\U_s^i).$$

\noindent Since $\U_s$ can be identified to an open subset of $\C^N$, and $\mu_s$ corresponds to a volume form proportional to the Lebesgue
measure, we have $\mu_s(\U_s^i)=0, \; \forall s \in \Ad, \; i\in\{1,\dots,N_1\}$. It follows immediately that $\mu_\Tr(\MET^{\mathrm{sp}})=0$.
\carre

\subsection{Proof of Theorem \ref{FVTh}}

From Corollary \ref{FVcor2}, and Proposition \ref{FVpr1}, to prove Theorem \ref{FVTh}, all we need is the following\\

\begin{proposition} \label{FVpr3}
Let $(\mathbf{A}_s;I;J)$, where $I=(i_1,\dots,i_N), J=(j_K,\dots,j_N)$, be an admissible triple. Let $\mathcal{F}_s$, and $\mu_s$ denote the
pull backs of the function $\mathcal{F}^\mathrm{et}$, and the volume form  $\mu_\Tr$ onto $\U_s$ by $\Phi_s$. Then we have:

$$\int_{\U_s(I;J)}\mathcal{F}_s d\mu_s < \infty.$$

\end{proposition}

\dem By the definition of primary system of indices,  we have a complex linear map

$$\begin{array}{crcc}
  \mathbf{B}_s :& \C^N & \lra & \ker\mathbf{A}_s \\
   & (z_1,\dots,z_N) & \longmapsto & (f_1(z_1,\dots,z_N),\dots,f_{N_1}(z_1,\dots,z_N)) \\
\end{array}$$

\noindent which is an isomorphism, where $f_{i_k}(z_1,\dots,z_N)=z_k$.  Consider a vector $(w_1,\dots,w_{N_1})$ in $\U_s$, let
$(\Sig,\hat{A},\xi)$ denote its image under $\Phi_s$. Let $\T, \hSig,\hT$ be as in the previous subsection. As usual, we denote the edges of
$\hT$ by $e_i, \; i=1,\dots,N_1$, so that $w_i$ is the complex number associated to $e_i$.\\

%

\noindent By definition, for any $k=K,\dots,N$, the complex numbers $w_{i_k}$ and $w_{j_k}$ correspond to two edges $e_{i_k}$, and $e_{j_k}$
which are contained in the boundary of a triangle $\Delta_k$ of $\hT$. With appropriate choices of orientation of $e_{i_k}$, and $e_{j_k}$, the
area of $\Delta_k$ is given by the function

$$\hat{\eta}_k=\frac{1}{2}(\mathrm{Re}(w_{i_k})\mathrm{Im}(w_{j_k})-\mathrm{Im}(w_{i_k})\mathrm{Re}(w_{j_k})).$$

\noindent Observe that the triangles $\Delta_k, \; k=K,\dots,N,$ are all distinct. Suppose on the contrary that there exist $k<k'$ such that
$e_{i_k},e_{j_k},e_{i_{k'}}$ are contained in the boundary of the same triangle. Excluding the cases $e_{i_{k'}}=e_{i_k}$, and
$e_{i_{k'}}=e_{j_k}$, we see that $e_{i_k},e_{j_k}$, and $e_{i_{k'}}$ are three sides of a triangle in $\hT$, which implies

$$w_{i_{k'}}=\pm w_{i_k}\pm  w_{j_k}.$$

\noindent Since $w_{j_k}$ is linearly dependent on $(w_{i_1},\dots,w_{i_{k-1}})$, it follows that $w_{i_{k'}}$ is linearly dependent on
$(w_{i_1},\dots,w_{i_k})$ as $(w_1,\dots,w_{N_1})$ varies in $\ker\As$, which is impossible since $(i_1,\dots,i_N)$ is a primary system of
indices for $\As$. As a consequence, we have

\begin{equation}\label{AreaIneq0}
\mathbf{Area}(\Sig) \geq \sum_{k=K}^N \hat{\eta}_k
\end{equation}

\noindent Let $\eta_k, \; k=K,\dots,N,$ denote the pull back of the function $\hat{\eta}_k$ by $\mathbf{B}_s$. It follows that
$\mathbf{B}_s^{-1}(\U_s(I;J))$ is contained in the following subset of $\C^{N}$

$$\mathcal{W}_s=\{(z_1,\dots,z_N)\in \C^N : \; |\mathrm{Re}(z_k)|\leq|\mathrm{Re}(f_{j_k}(z_1,\dots,z_N))|,
\; \eta_k>0, \;  k=K,\dots,N\}.$$

\noindent Let $\mathcal{G}_s$ denote the pull back of $\mathcal{F}_s$ by $\mathbf{B}_s$, since the volume form $\mathbf{B}_s^*\mu_s$ equals
$\kappa\lambda_{2N}$, where $\lambda_{2N}$ is the Lebesgue measure of $\C^N$, and $\kappa$ is a constant, all we need to show is the following

\begin{equation}\label{InEq1}
\int_{\mathcal{W}_s}\mathcal{G}_s d\lambda_{2N} <\infty
\end{equation}


To simplify the notations, for $k=1,\dots,N$, set $x_k=\mathrm{Re}(z_k), \; y_k=\mathrm{Im}(z_k)$. For $k=K,\dots,N$, we write $f_k$ instead of
$f_{j_k}$, and set $a_k=\mathrm{Re}(f_k),\; b_k=\mathrm{Im}(f_k)$. Recall that, by definition, $f_k$ depends only on $(z_1,\dots,z_{k-1})$, and
the coefficients of $z_K,\dots,z_{k-1}$ in $f_k$ are all real. Thus, we deduce that $a_k$ is a function of
$(z_1,\dots,z_{K-1},x_K,\dots,x_{k-1})$, and $b_{k}$ is a function of $(z_1,\dots,z_{K-1},y_K,\dots,y_{k-1})$, for $k=K,\dots,N$. With these
notations, we have

\begin{equation}\label{TriArea}
\eta_k = \frac{1}{2}(x_kb_k-y_ka_k), \; k=K,\dots,N.
\end{equation}

\begin{equation}\label{MajCond}
|x_k|  \leq   |a_k|, \;  k=K,\dots,N.
\end{equation}

\begin{equation}\label{AreaIneq}
\mathbf{Area}(\Sig)\geq\sum_{k=K}^N\eta_k.
\end{equation}

Recall that, by definition of admissible triple, the complex numbers $z_1,\dots, z_{K-1}$ correspond to some edges of $\hT$ in the boundary of
$\hSig$, or equivalently to some edges of the forest $\hat{A}$.  Therefore, we have

\begin{equation}\label{BdyIneq}
\ell^2(\hat{A})\geq  \sum_{k=1}^{K-1} |z_k|^2
\end{equation}

\noindent Consequently, we have

\begin{equation}\label{InEq2}
\mathcal{G}_s\leq \exp(-\sum_{k=1}^{K-1}|z_k|^2-\sum_{k=K}^N\eta_k)
\end{equation}

\noindent Therefore, to prove (\ref{InEq1}), it suffices to show

\begin{lemma}\label{FVlm2}
\begin{equation}\label{InEq3}
\mathcal{I}=\int_{\mathcal{W}_s} \exp(-\sum_{k=1}^{K-1}|z_k|^2-\sum_{k=K}^N\eta_k) d\lambda_{2N}<\infty
\end{equation}
\end{lemma}

\dem Fix $(z_1,\dots,z_{K-1})\in \C^{K-1}$ and $(x_{K},\dots,x_{N})\in \R^{N-K+1}$, and set

$$\begin{array}{cl}
\mathcal{W}_s((z_1,\dots,z_{K-1});(x_{K},\dots,x_{N}))= & \{(y_{K},\dots,y_N) \in \R^{N-K+1} \text{ such that}\\
 & (z_1,\dots,z_{K-1},x_K+\imath y_{K},\dots,x_N+\imath y_N) \in \mathcal{W}_s\}\\
\end{array}$$

\noindent Consider the following integral

$$\mathcal{I}((z_1,\dots,z_{K-1});(x_{K},\dots,x_{N}))= \int_{{\mathcal{W}_s}((z_1,\dots,z_{K-1});(x_{K},\dots,x_{N}))}\exp(-\sum_{k=
K}^N\eta_k)dy_{K}\dots dy_N.$$

\noindent Making the change of variables $(y_{K},\dots,y_{N})\mapsto(\eta_{K},\dots,\eta_N)$,  using (\ref{TriArea}), and the fact that, with
$(z_1,\dots,z_{K-1},x_K,\dots,x_N)$ fixed, $a_k$ is constant, and $b_{k}$ is an affine function of $(y_K,\dots,y_{k-1})$, for any $k=K,\dots,N$,
we have:

$$d\eta_{K}\dots d\eta_N= \frac{|a_{K}\dots a_{N}|}{2^{N-K+1}}dy_{K}\dots dy_N.$$

\noindent Since for  $k=K,\dots,N$, the function $\eta_k$ is positive on $\mathcal{W}_s$, it follows

\begin{eqnarray*}
  \mathcal{I}((z_1,\dots,z_{K-1});(x_{K},\dots,x_{N})) & \leq & \frac{2^{N-K+1}}{|a_{K}\dots a_{N}|}\int_0^{+\infty}
e^{-\eta_{K}}d\eta_{K}\dots\int_0^{+\infty} e^{-\eta_N}d\eta_N \\
   &\leq & \frac{2^{N-K+1}}{|a_{K}\dots a_{N}|}.\\
\end{eqnarray*}

\noindent Now, set

$$\mathcal{W}^*_s=\{((z_1,\dots,z_{K-1});(x_{K},\dots,x_N))\in \C^{K-1}\times\R^{N-K+1} \; : \; |x_k|\leq |a_k|,\; k=K,\dots,N\}.$$

\noindent We have

\begin{eqnarray*}
\mathcal{I} &=&\int_{\mathcal{W}^*_s}\exp(-\sum_{k=1}^{K-1}|z_k|^2)\mathcal{I}((z_1,\dots,z_{K-1});(x_{K},\dots,x_{N}))dx_1dy_1\dots
dx_{K-1}dy_{K-1}dx_{K}\dots dx_N,\\
&\leq&\int_{\mathcal{W}^*_s}\exp(-\sum_{k=1}^{K-1}|z_k|^2) \frac{2^{N-K+1}}{|a_{K}\dots a_{N}|}dx_1dy_1\dots
dx_{K-1}dy_{K-1}dx_{K}\dots dx_N,\\
&\leq&\int_{\C^{K-1}}\exp(-\sum_{k=1}^{K-1}|z_k|^2)[\int_{-|a_{K}|}^{|a_{K}|}[\dots[\int_{-|a_N|}^{|a_N|}
\frac{2^{N-K+1}}{|a_{K}\dots a_{N}|} dx_N]\dots]dx_{K}]dx_1dy_1\dots dx_{K-1}dy_{K-1}.\\
\end{eqnarray*}

\noindent Using the fact that $a_k$ does not depend on $x_k,\dots,x_N$ for $k=K,\dots,N$, we see that

$$\int_{-|a_{K}|}^{|a_{K}|}[\dots[\int_{-|a_N|}^{|a_N|}\frac{2^{N-K+1}}{|a_{K}\dots a_{N}|} dx_N]\dots]dx_{K}=4^{N-K+1}.$$

\noindent Hence,

$$\mathcal{I} \leq 4^{N-K+1}\int_{\C^{K-1}} e^{-(|z_1|^2+\dots+|z_{K-1}|^2)}dx_1dy_1\dots dx_{K-1}dy_{K-1} <\infty.$$

\noindent The lemma is then proved.\carre


\noindent The proofs of Proposition \ref{FVpr3}, and of Theorem \ref{FVTh} are now complete. \carre

\section{Finiteness of the volume of moduli spaces of translation surfaces}\label{prAPrf}


In this section, we prove Theorem \ref{FVprA} using Theorem \ref{FVTh}. Recall that $\Hg$ is the moduli space of triples
$(\Sig,\{x_1,\dots,x_n\},\xi)$, where $\Sig$ is a translation surface, $\{x_1,\dots,x_n\}$ is the set of marked singularities of $\Sig$ with
cone angles $\{(k_1+1)2\pi,\dots,(k_n+1)2\pi\}$ respectively, and $\xi$ is a unitary parallel vector filed on the complement of the set
$\{x_1,\dots,x_n\}$. An element of $\Hg$ can be identified to a pair $(M,\omega)$, where $M$ is a connected, closed Riemann surface, and
$\omega$ is a holomorphic $1$-form on $M$ having exactly $n$ zeros with orders $k_1,\dots,k_n$. Using this identification, one can define a
local chart for $\Hg$ in a neighborhood of a point $(M,\omega)$ as follows: let $p_1,\dots,p_n$ denote the zeros of $\omega$, and let
$(\gamma_1,\dots,\gamma_{2g+n-1})$ be a basis of $H_1(M,\{p_1,\dots,p_n\},\Z)$. There exists a neighborhood $\U$ of $(M,\omega)$ in $\Hg$ such
that, for any $(M',\omega')$ in $\U$, $(\gamma_1,\dots,\gamma_{2g+n-1})$ gives rise to a basis $(\gamma'_1,\dots,\gamma'_{2g+n-1})$ of
$H_1(M',\{p'_1,\dots,p'_n\},\Z)$, where $p'_1,\dots,p'_n$ are the zeros of $\omega'$. It follows that the map

$$\begin{array}{cccc}
  \Phi: & \U & \lra & \C^{2g+n-1} \\
    & (M',\omega') & \longmapsto & (\int_{\gamma'_1}\omega',\dots,\int_{\gamma'_{2g+n-1}}\omega') \\
\end{array}$$

\noindent is a local chart. Let $\mu_0$ denote the pull-back of the Lebesgue measure of $\C^{2g+n-1}$ by $\Phi$, then $\mu_0$ is a well-defined
volume form on $\Hg$.\\

Since $\Hg$ is a special case of flat surface with erasing forest, where all the trees in the forest are points, on $\Hg$, we also have a volume
form $\mu_\Tr$.  It turns out that(cf. \cite{Ng1}), on each connected component of $\Hg$, we have $\mu_\Tr=\lambda \mu_0$, where $\lbd$ is a
constant.\\

\subsection{Translation surface with a marked geodesic segment}

To prove Theorem \ref{FVprA}, we first consider a space $\MET$ where all the trees in $\Ah$ but one are points, and the remaining one is a
segment, together with a projection $\DS{\varrho: \MET \lra \Hg }$ which is (locally) a fiber bundle.  We then use Theorem \ref{FVTh}, and the
fact that the integral of $\mathcal{F}^\et$ on the fibers of $\varrho$ is constant to conclude.\\

Set $\alpha_i=2(k_i+1),\; i=1,\dots,n$. Let $\A_1$ be a topological tree isomorphic to a segment, and for $i=2,\dots,n$, let $\A_i$ be just a
point. Let $\underline{\alpha}$ denote the vector $(2\pi, \alpha_1,\dots,\alpha_n)$, and $\Ah$ denote the family $\{\A_1,\dots,\A_n\}$. Consider
the space $\MET$ with the previous data. In this case, $\MET$ is the moduli space of triples $(\Sig,\{I(x_1,x),x_2,\dots,x_n\},\xi)$, where

\begin{itemize}
\item[.]  $\Sig$ is a translation surface,

\item[.] $\{x_1,\dots,x_n\}$ is the set of singularities of $\Sig$ with cone angles $\{\alpha_1,\dots,\alpha_n\}$ respectively,

\item[.] $x$ is a regular point of $\Sig$,

\item[.] $I(x_1,x)$ is a geodesic segment joining the singular point $x_1$ to $x$,

\item[.] and $\xi$ is a unitary parallel vector field on the complement of $I(x_1,x)\cup\{x_2,\dots,x_n\}$.

\end{itemize}

By definition, we have a natural projection $\varrho$ from $\MET$ to $\Hg$ consisting of forgetting the segment $I(x_1,x)$, that is

$$ \varrho: (\Sig,\{I(x_1,x),x_2,\dots,x_n\},\xi) \longmapsto (\Sig,\{x_1,\dots,x_n\},\xi).$$

Let $\DS{N=\dim_\C \Hg}$, clearly, $\dim_\C\MET=N+1$. Let $\hat{\mu}_\Tr$, and $\mu_\Tr$ denote the volume forms on $\MET$ and $\Hg$ arising
from admissible triangulations respectively.\\

Let $\Phi$ denote the period mapping defining a local chart of $\Hg$ in an open set $\U$. We can then define some local charts $\hat{\Phi}$ for
$\MET$ whose domains cover $\varrho^{-1}(\U)$ as follows: first, we identify any  $(\Sig,\{x_1,\dots,x_n\},\xi)$ in $\U$ to a pair $(M,\omega)$,
if $\Phi(\Sig,\{x_1,\dots,x_n\},\xi)=(z_1,\dots,z_N)$, then $\hat{\Phi}(\Sig,\{I(x_1,x),x_2,\dots,x_n\},\xi)= (z_1,\dots,z_{N},z_{N+1})$, where

$$z_{N+1}=\int_{x_1}^x\omega,$$

\noindent and the integral is taken along $I(x_1,x)$. In the local charts $\hat{\Phi}$, and $\Phi$, the map $\varrho$ can be written as

\begin{equation*}
\varrho(z_1,\dots,z_{N+1})=(z_1,\dots,z_N).
\end{equation*}

Let $\lbd_{2N}$, and $\lbd_{2(N+1)}$ denote the Lebesgue measures of $\C^N$, and $\C^{N+1}$. Up to some multiplicative constants depending on
the connected component of $(\Sig,\{x_1,\dots,x_n\},\xi)$ in $\Hg$, we can write

\begin{equation}\label{prAeq0}
\mu_\Tr=\Phi^*\lbd_{2N} \text{ and } \hat{\mu}_\Tr=\hat{\Phi}^*\lbd_{2(N+1)}.
\end{equation}

\subsection{ Proof of Theorem \ref{FVprA}}

Consider a point $(\Sig,\{x_1,\dots,x_n\},\xi)$ in $\Hg$. Fix a unitary tangent vector $v_1\in T_{x_1}\Sig$, we can then identify the set of
unitary tangent vectors  of $T_{x_1}\Sig$ to $\R/(\alpha_1\Z)$. Any geodesic segment in $\Sig$ which contains $x_1$ as an endpoint is uniquely
determined by its tangent vector at $x_1$, and its length. Hence, we have an injective map:

$$ \varphi :\varrho^{-1}\{(\Sig,\{x_1,\dots,x_n\},\xi)\} \lra (\R/(\alpha_1\Z))\times\R^+.$$

\noindent Let $\U$ be a neighborhood of $(\Sig,\{x_1,\dots,x_n\},\xi)$ in $\Hg$ on which a local chart can be defined by some period mapping
$\Phi$. For each point $(\Sig',\{x'_1,\dots,x'_n\},\xi')$ in $\U$, we choose a tangent vector $v'_1$ in $T_{x'_1}\Sig'$ to be the reference
vector, we can assume that $v'_1$ varies continuously as $(\Sig',\{x'_1,\dots,x'_n\},\xi')$ varies in $\U$ so that the map $\varphi$ can be
extended into an injective, continuous  map:

$$\varphi:\varrho^{-1}(\U)\lra \U\times(\R/\alpha_1\Z)\times\R^+.$$

\noindent Using the local charts $\hat{\Phi}$ on $\varrho^{-1}(\U)$, we can write

\begin{equation}\label{prAeq0b}
\varphi(z_1,\dots,z_{N+1})=((z_1,\dots,z_N),\mathrm{arg}(z_{N+1})+c, |z_{N+1}|), \;  \text{ where } c \text{ is some constant}
\end{equation}

\noindent Let $d\theta$, and $dr$ denote the standard measures on $\R/(\alpha_1\Z)$, and $\R^+$ respectively. From (\ref{prAeq0}), and
(\ref{prAeq0b}), we have

$$\varphi_*d\hat{\mu}_\Tr=rd\mu_\Tr d\theta dr.$$

\noindent Consequently,

\begin{equation}\label{prAeq4}
\int_{\varrho^{-1}(\U)}e^{-\mathbf{Area}(\Sig)-\ell^2(I)}d\hat{\mu}_\Tr=\int_{\varphi(\varrho^{-1}(\U))}e^{-\mathbf{Area}(\Sig)-r^2}r d\mu_\Tr
d\theta dr.
\end{equation}

\noindent By a well known result (for example, see \cite{MasTab}, Theorem 1.8), we know that, on a translation surface, there exists a countable
subset $\Theta$ of $\R/\alpha_1\Z$ such that if $\theta$ is not in $\Theta$, then the geodesic ray starting from $x_1$ in the direction $\theta$
can be extended infinitely. It follows immediately that $\varphi(\varrho^{-1}(\U))$ is an open dense subset, hence of full measure, of
$\U\times(\R/\alpha_1\Z)\times\R^+$. Therefore, we have

\begin{eqnarray*}
\int_{\varphi(\varrho^{-1}(\U))}e^{-\mathbf{Area}(\Sig)-r^2}rd\mu_\Tr d\theta dr &=&
\int_{\U\times(\R/\alpha_1\Z)\times\R^+} e^{-\mathbf{Area}(\Sig)-r^2}rd\mu_\Tr d\theta dr,\\
 &=&\int_0^{+\infty}e^{-r^2}rdr\int_0^{\alpha_1}d\theta\int_\U e^{-\mathbf{Area}(\Sig)}d\mu_\Tr,\\
 &=&\frac{\alpha_1}{2}\int_\U e^{-\mathbf{Area}(\Sig)}d\mu_\Tr.\\
\end{eqnarray*}

\noindent It follows from (\ref{prAeq4}) that

\begin{equation}\label{prAeq5}
\int_{\varrho^{-1}(\U)}e^{-\mathbf{Area}(\Sig)-\ell^2(I)} d\hat{\mu}_\Tr=\frac{\alpha_1}{2}\int_{\U}e^{-\mathbf{Area}(\Sig)}d\mu_\Tr
\end{equation}

\noindent Since (\ref{prAeq5}) is true for any small neighborhood $\U$ in $\Hg$, we deduce that

\begin{equation*}
\int_{\MET}e^{-\mathbf{Area}(\Sig)- \ell^2(I)}d\hat{\mu}_\Tr=\frac{\alpha_1}{2}\int_{\Hg}e^{-\mathbf{Area}(\Sig)}d\mu_\Tr.
\end{equation*}

\noindent By Theorem \ref{FVTh}, we know that

\begin{equation*}
\int_{\MET} e^{-\mathbf{Area}(\Sig)-\ell^2(I)}d\hat{\mu}_\mathrm{Tr} < \infty.
\end{equation*}

\noindent Therefore,

\begin{equation*}
\int_{\Hg}e^{-\mathbf{Area}(\Sig)}d\mu_\Tr <\infty,
\end{equation*}

\noindent and Theorem \ref{FVprA} is then proved. \carre

\section{Finiteness  of $\mu^1_\Tr(\MSSi)$}\label{prBPrfSect}

In this section, we are interested in the moduli space of flat surfaces of genus zero with prescribed cone angles. Let $\MSSv$ denote the moduli
space of flat surfaces having $n$ singularities, which are numbered, with cone angles given by $\underline{\alpha}=(\alpha_1,\dots,\alpha_n)$.
Recall that we have a volume form $\mu_\mathrm{Tr}$ on the space $\MSS=\MSSv\times\S^1$, which is defined by identifying locally $\MSS$ to
$\MET$, with some appropriate choice of $\Ah$.  Let $\MSSvi$ de the set of surfaces having unit area in $\MSSv$, and $\MSSi$ denote the product
space $\MSSvi\times \S^1$. The space $\MSSvi$ can be considered as the moduli space of the configurations of $n$ marked points on the sphere
$\S^2$ up to Möbius transformations.\\

\noindent The volume form $\mu_\mathrm{Tr}$ induces naturally a volume form $\mu^1_\mathrm{Tr}$ on the space $\MSSi=\mathbf{Area}^{-1}(\{1\})$,
and hence,  a volume form $\hat{\mu}^1_\mathrm{Tr}$ on $\MSSvi$ by pushing forward. As we have seen in the introduction, Theorem \ref{FVprB} is
equivalent to the finiteness of  $\mu^1_\Tr(\MSSi)$, and of $\hat{\mu}^1_\Tr\MSSvi$. Our aim in this section is to prove Theorem \ref{FVprB}
using Theorem \ref{FVTh}.\\



\subsection{The function $\delta$}

Let $\Sig$ be an element of $\MSSv$. Let $x_1,\dots,x_n$ denote the singular points of $\Sig$ so that the cone angle at $x_i$ is $\alpha_i$. Let
$\di$ denote the distance induced by the flat metric on $\Sig$. For any subset $I$ of $\{1,\dots,n\}$, let $\diam_I(\Sig)$ denote the diameter
of the set $\{x_i,\; i\in I\}$.  We define

\begin{equation*}
\delta_I(\Sig)= \min\{\di(x_i,x_j): \; i \in I, \; j \not\in I\},\\
\end{equation*}

\noindent and

\begin{equation*}
\delta^+_I(\Sig)=\left\{ \begin{array}{l}
\delta_I(\Sig)   \hbox{ \hspace{0.5cm} if  $\delta_I(\Sig) \geq 3 \diam_I(\Sig)$ ,}\\
0   \hbox{   \hspace{0.5cm} otherwise.}\\
\end{array} \right.
\end{equation*}

\noindent A subset $I$ of $\{1,\dots,n\}$ is called {\em essential} if we have

$$\sum_{i\in I}\alpha_i\not\in 2\pi\N.$$

\noindent We define a function $\delta$ on the space $\MSSv$ as follows

$$\text{ for every } \Sig \in \MSSv, \; \delta(\Sig)=\max\{\delta^+_I(\Sig) : \; I \subset \{1,\dots,n\}, \;  I \text{ is essential }\}.$$

\rem The function $\delta$ is always positive, since when $I=\{i\}$, we have

$$\delta^+_{\{i\}}(\Sig)=\min\{\di(x_i,x_j),\; j \neq i\} >0,$$

\noindent and there always exists $i\in \{1,\dots,n\}$ such that $\alpha_i \not\in 2\pi\N$, which means that $\{i\}$ is essential. To simplify
the notations, we also denote by $\delta$ the composition of $\delta$ with the natural projection from $\MSS=\MSSv\times\S^1$ onto $\MSSv$.\\

\subsection{Good tree and good forest}\label{GdForSect}

Fix a surface  $\Sig$ in $\MSSv$, and let $x_1,\dots,x_n$ denote the singular points of $\Sig$ so that the cone angle at $x_i$ is $\alpha_i$.
Let $\mathrm{V}$ denote the set $\{x_1,\dots,x_n\}$, and set $\delta=\delta(\Sig)$. For any geodesic tree $A$ on $\Sig$, we denote by $\ver(A)$
the vertex set of $A$, $\max(A)$ the length of the longest edge of $A$, and by $\Ra(A)$ the distance from $\ver(A)$ to the set $\mathrm{V}\setm
\ver(A)$.\\

\begin{definition}\label{DefGdTree}
Let $A$ be a geodesic tree in $\Sig$ whose vertex set is a subset of $\mathrm{V}$. Let $k$ be the number of edges of $A$. The tree $A$ is said
to be {\em good}, if either $A$ is a singular point with cone angle in $2\pi\N$, or $k\geq 1$ and we have

\begin{itemize}
\item[$\bullet$] $\max(A)\leq 4^{k-1}\delta$,

\item[$\bullet$] $\diam(\ver(A))\leq 4^{k-1}\delta$,

\item[$\bullet$] The set of indices corresponding to the vertex set of $A$ is non essential, that is the sum of all cone angles at the vertices
of $A$ belongs to $2\pi\N$.

\item[$\bullet$] Either $\ver(A)=\mathrm{V}$, or $\Ra(A)\geq 3.4^{k-1}\delta$.

\end{itemize}

\noindent A union of disjoint good trees such that the union of the vertex sets is $\mathrm{V}$ is called a  {\em good forest}.\\

\end{definition}

We have

\begin{lemma}\label{FVSSlm1}
There always exists a good forest in $\Sig$.\\
\end{lemma}

\noindent The proof of this lemma is given in Appendices, Section \ref{FVSSlm1Sect}.

\begin{corollary}\label{FVSScor1}
There exists a constant $\kappa$ depending only on $n$ such that for any $\Sig$ in $\MSSv$, there exists an erasing forest $\hat{A}$ in $\Sig$
which verifies

$$\ell(\hat{A})\leq  \kappa\delta.$$

\end{corollary}

\dem By Lemma \ref{FVSSlm1}, we know that there exists a good forest $\hat{A}=\sqcup_{j=1}^m A_j$ in $\Sig$. By definition, for every $j \in
\{1,\dots,n\}$, the sum of the cone angles at the vertices of $A_j$ belongs to $2\pi\N$, therefore, $\hat{A}$ is an erasing forest. Since every
tree $A_j$ in $\hat{A}$ is good, we have

$$\ell(A_j)\leq k_j4^{k_j-1}\delta,$$

\noindent where $k_j$ is the number of edges of $A_j$. Observe that $k_1+\dots+k_m=n-m\leq n-1$. Therefore, we have

$$\ell(\hat{A})=\sum_{j=1}^m\ell(A_j)\leq (n-1)4^{n-1}\delta,$$

\noindent and the corollary follows.\carre

\subsection{Proof of Theorem \ref{FVprB}}\label{prfFVprB}

Theorem \ref{FVprB} is a consequence of two following propositions:

\begin{proposition}\label{FVSSpr1}
We have

\begin{equation*}
\int_{\MSS}\exp(-\mathbf{Area}-\delta^2)d\mu_\mathrm{Tr} <\infty.
\end{equation*}

\end{proposition}

\noindent and

\begin{proposition}\label{FVSSpr2}
There exists a constant $C(\underline{\alpha})$ depending on $\underline{\alpha}$ such that for any surface $\Sig$ in $\MSSv$ we have

$$\delta^2(\Sig)< C(\underline{\alpha})\mathbf{Area}(\Sig).$$
\end{proposition}

\noindent The proof of Proposition \ref{FVSSpr2} is rather straight forward but quite lengthy, it will be given in Appendices, Section
\ref{FVSSpr2Sect}. Here below, we  give the proof of Proposition \ref{FVSSpr1} using Corollary \ref{FVSScor1}.\\



\dem (of Proposition \ref{FVSSpr1}) Let $\A_\mathrm{ad}(\underline{\alpha})$ denote the set of all families $\Ah=\{\A_1,\dots,\A_m\} \; (0<m<n)$
of topological trees, whose vertices are labelled by $\{1,\dots,n\}$, up to isomorphism, verifying the following condition: if $I_j, \;
j=1,\dots,m,$ is the subset of $\{1,\dots,n\}$ in bijection with the vertices of the tree $\A_j$, then

\begin{equation*}
\sum_{i\in I_j}\alpha_i \in 2\pi\N.
\end{equation*}

\noindent For each $\Ah=\{\A_1,\dots,\A_m\} \in \A_\mathrm{ad}(\underline{\alpha})$, let $\U_{\Ah}$ denote the subset of $\MET$ consisting of
all triples $(\Sig,\hat{A},\xi)$ satisfying the following condition:

\begin{equation*}
\ell(\hat{A})\leq \kappa\delta(\Sig),\\
\end{equation*}

\noindent where  $\kappa$ is the constant in Corollary \ref{FVSScor1}. Let $\rho_{\Ah}$ denote the map from $\MET$ onto $\MSSv$, which
associates to every triple $(\Sig,\hat{A},\xi)$ the surface $\Sig$. From Corollary \ref{FVSScor1}, we know that the family

$$\{\V_{\Ah}=\rho_{\Ah}(\U_{\Ah}): \; \Ah \in \A_\mathrm{ad}(\underline{\alpha})\}$$

\noindent covers the space $\MSSv$. Let $\rho_1$ be the natural projection from $\MSS$ onto $\MSSv$, it follows that the family

$$\{\rho_1^{-1}(\V_{\Ah}): \; \Ah \in \A_\mathrm{ad}(\underline{\alpha})\}$$

\noindent covers the space $\MSS$. Since the set $\A_\mathrm{ad}(\underline{\alpha})$ is finite, it is enough to show that, for every $\Ah$ in
$\A_\mathrm{ad}(\underline{\alpha})$, we have

\begin{equation}\label{FVSSeq1}
\int_{\rho_1^{-1}(\V_{\Ah})}\exp(-\mathbf{Area}-\delta^2)d\mu_\mathrm{Tr}<\infty.
\end{equation}

\noindent Since the space $\MSS$ can be locally identified to $\MET$, we have

\begin{equation*}
\int_{\rho_1^{-1}(\V_{\Ah})}\exp(-\mathbf{Area}-\delta^2)d\mu_\mathrm{Tr}=\int_{\U_{\Ah}}\exp(-\mathbf{Area}-\delta^2)d\mu_\mathrm{Tr}
\end{equation*}

\noindent By definition, for every $(\Sig,\hat{A},\xi)$ in $\U_{\Ah}$, we have $\ell(\hat{A})\leq \kappa\delta(\Sig)$. It follows

\begin{equation}\label{FVSSeq2}
\int_{\U_{\Ah}}\exp(-\mathbf{Area}-\delta^2)d\mu_\mathrm{Tr}\leq \int_{\U_{\Ah}}\exp(-\mathbf{Area}-\frac{1}{\kappa^2}\ell^2)d\mu_\mathrm{Tr}
\end{equation}

\noindent By Theorem \ref{FVTh}, we know that the right hand side of (\ref{FVSSeq2}) is finite. Consequently, (\ref{FVSSeq1}) is true, and the
proposition follows.\carre

Proposition \ref{FVSSpr1}, and Proposition \ref{FVSSpr2} imply that

\begin{equation}\label{FVSSeq3}
\int_{\MSS} \exp(-(1+C(\underline{\alpha}))\mathbf{Area}(.))d\mu_\Tr < \infty
\end{equation}

\noindent which is equivalent to

\begin{equation}\label{FVSSeq4}
\int_{\MSS} \exp(-\mathbf{Area}(.))d\mu_\Tr < \infty
\end{equation}

\noindent since both (\ref{FVSSeq3}), and (\ref{FVSSeq4}) are equivalent to the fact that the volume of $\MSSi$ is finite. The proof of Theorem
\ref{FVprB} is now complete. \carre


\begin{appendices}

\section{Existence of good forest}\label{FVSSlm1Sect}

\subsection{ Existence of good tree}

Let $\Sig, x_1,\dots,x_n, \mathrm{V}, \delta$ be as in Section \ref{GdForSect}. Let $\di$ denote the distance induced by the metric of $\Sig$.
Let us start by proving the following

\begin{lemma}\label{FVSSlm2}
For any $\Sig$ in $\MSSv$, there always exists a good tree on $\Sig$.
\end{lemma}

\dem First, let $e$ be a geodesic segment which realizes the distance

$$\min\{\di(x_i,x_j),\; \alpha_i\not\in 2\pi\N  \text{ and } i\neq j \}.$$

\noindent  By definition, we have

$$ \leng(e)=\min\{\delta^+_{\{i\}}(\Sig),\; \alpha_i\not\in 2\pi\N\} \leq \delta.$$

\noindent Let $A^1$ denote the tree which contains only the segment $e$. By assumption, we have

$$\max(A^1)=\diam(\ver(A^1))=\leng(e_1)\leq \delta.$$

\noindent Consider the following procedure, which will be called the {\em vertex adding procedure}: suppose that we already have a geodesic tree
$A^k, \; k\geq 1,$ connecting $k+1$ points in $\{x_1,\dots,x_n\}$ verifying the following condition:

$$ (*)\left\{ \begin{array}{ccr}
\max(A^k) & \leq & 4^{k-1}\delta,\\
\diam(\ver(A^k)) & \leq & 4^{k-1}\delta.\\
\end{array} \right.$$

\noindent Let  $I$ be the subset of $\{1,\dots,n\}$ corresponding to the vertex set of $A^k$. We have two cases:

\begin{itemize}
\item[-]\underline{Case 1:} $I$ is essential. In this case, let $e_{k+1}$ be a segment realizing the distance $\delta_I(\Sig)$, and let $x_{j}$
be the endpoint of $e_{k+1}$ which does not belong to $\ver(A^k)$. By definition, we have either

\begin{itemize}
\item[.]  $\leng(e_{k+1})\leq 3\diam(\ver(A^k))$  or,

\item[.]  $\leng(e_{k+1})\leq \delta$.
\end{itemize}

\noindent Since $\diam(\ver(A^k))\leq 4^{k-1}\delta$, it follows that $\DS{\leng(e_{k+1})\leq 3.4^{k-1}\delta}$,  in both cases

Slit open the surface $\Sig$ along the tree $A^k$, and let $\Sig'$ denote the new surface. The vertex set $\ver(A^k)$ of $A^k$ gives rise to a
finite subset $V^k$ of the boundary of $\Sig'$. Let us prove that the distance in $\Sig'$ from $x_j$ to $V^k$ is at most $4^k\delta$.\\

Consider $e_{k+1}$ as a ray exiting from $x_{j}$, and let $y$ be the first intersection point between $e_{k+1}$ and the tree $A^k$. Since we
have $\max(A^k)\leq 4^{k-1}\delta$, there exists a path in $\Sig$ joining $x_j$ to an endpoint of the edge containing $y$ without crossing any
edge of $A^k$, whose length is at most

$$3.4^{k-1}\delta+4^{k-1}\delta=4^k\delta.$$

Because this path does not cross any edge of the tree $A^k$, it represents a path on $\Sig'$ joining $x_j$ to a point in $V^k$. Thus, we deduce
that the distance between $x_j$ and $V^k$ in $\Sig'$ is at most $4^k\delta$.\\

\noindent The path realizing the distance from $x_j$ to $V^k$ in $\Sig'$ corresponds to a path $a$ in $\Sig$ which is piecewise geodesic with
endpoints in $\mathrm{V}$, joining $x_j$ to a vertex of the tree $A^k$. Note that we have

$$\leng(a)\leq 4^{k}\delta.$$

Adding $ a$ to $A^k$, we get a new tree which contains $k+r$ edges, and will be denoted by $A^{k+r}$, where $r$ is the number of geodesic
segments with endpoints in $\mathrm{V}$ contained in $a$. Let us prove that  this new tree also verifies the condition $(*)$.

\begin{itemize}

\item[$\bullet$] If $r=1$, then $\ver(A^{k+1})=\ver(A^k)\cup\{x_j\}$. Since $\diam(A^k)\leq 4^{k-1}\delta$, and the distance from $x_j$ to
$\ver(A^k)$ is at most $3.4^{k-1}\delta$, we deduce that

$$\diam(\ver(A^{k+1}))\leq 4^{k-1}\delta+3.4^{k-1}\delta=4^k\delta.$$

By assumption, we know that $\max(A^k)\leq 4^{k-1}\delta$, and we have proved that the length of the added edge is at most $4^k\delta$, hence,
we have $\max(A^{k+1}) \leq 4^k\delta$.\\

\item[$\bullet$] If $r>1$, it means that the path $a$ contains some points of $\mathrm{V}$ in its interior. The distance from these points to
the set $\ver(A^k)$ is bounded by the length of $a$ which is at most $4^k\delta$. Hence, the diameter of the set $\ver(A^{k+r})$ is at most

$$4^{k-1}\delta+4^k\delta\leq 4^{k+r-1}\delta.$$

As for $\max(A^{k+r})$, we have

$$\max(A^{k+r})\leq\max\{\max(A^k),\leng(a)\}\leq 4^k\delta.$$

\end{itemize}

We can now restart the procedure with $A^{k+r}$ in the place of $A^k$.\\

\item[-] \underline{Case 2:} $I$ is non-essential. In this case, if $\ver(A^k)=\mathrm{V}$, or $\Ra(\ver(A^k))\geq 3.4^{k-1}\delta$, then the
procedure stops since we already get a good tree. Otherwise, there exist $x_i$ in $\ver(A^k)$, $x_j$ in $\mathrm{V}\setm \ver(A^k)$, and a
geodesic segment $e$ joining $x_i$ to $x_j$ with

$$\leng(e) < 3.4^{k-1}\delta.$$

Using the same arguments as in Case 1, we can add to $A^k$ some edges so that the new tree also verifies the condition $(*)$, and repeat the
procedure.\\
\end{itemize}

\noindent Since we only have finitely many singular points in $\Sig$, the vertex adding procedure must stop, and we obtain a good tree. \carre


\subsection{ Proof of Lemma \ref{FVSSlm1}}

By Lemma \ref{FVSSlm2}, we know that there exists a good tree $A_1$ in $\Sig$. If $\ver(A_1)=\mathrm{V}$, or every point in the set
$\mathrm{V}\setm\ver(A_1)$ has cone angle in $2\pi\N$, then we are done. Otherwise, there exists a point $x_i$ in $\mathrm{V}\setm\ver(A_1)$,
with cone angle not in the set $2\pi\N$. In this case, we would like to construct a good tree $A_2$ containing $x_j$ by the vertex adding
procedure. However, this procedure can not be carried out straightly because of the presence of the tree $A_1$. Namely, it may happen that we
have $\Ra(\ver(A_2))\leq 3.4^{k_2-1}\delta$, where $k_2$ is the number of edges of $A_2$, but the segment realizing the distance $
\di(\ver(A_2), \mathrm{V}\setm \ver(A_2))$ intersects the tree $A_1$.\\

\noindent To fix this problem, let us consider the following procedure, which will be called the {\em tree joining procedure}: let
$A_1,\dots,A_{l}$ be a family of disjoint geodesic trees whose vertex sets are contained in $\mathrm{V}$. Let $k_1,\dots,k_l,\; k_i>0,$ be the
numbers of edges of $A_1,\dots,A_{l}$ respectively. Assume that the family $\{A_1,\dots,A_{l}\}$ verifies the following properties:

$$ (**)\left\{%
\begin{array}{ll}
    a) & \hbox{ $A_1,\dots,A_{l-1}$  are good trees,} \\
    b) & \hbox{ $A_l$  satisfies the condition $(*)$,} \\
    c) & \hbox{ $\di(A_l,\sqcup_{j=1}^{l-1}A_j)\leq 3.4^{k_l-1}\delta$.} \\
\end{array}%
\right.$$

\noindent Let  $s$ be a path of length at most $3.4^{k_l-1}\delta$ joining a point of $A_l$ to a point of $\sqcup_{j=1}^{l-1} A_j$. Without loss
of generality, we can assume that $s$ joins a point in $A_l$ to a point in $A_{l-1}$. Since both $A_{l-1}$ and $A_l$ verify the condition $(*)$,
in particular, we have

$$\max(A_l)\leq 4^{k_l-1}\delta, \text{ and } \max(A_{l-1})\leq 4^{k_{l-1}-1}\delta.$$

\noindent It follows that there exists a path $c$ joining a vertex of $A_{l-1}$ to a vertex of $A_l$ without crossing any edge of the family
$\{A_1,\dots,A_l\}$ such that

$$\leng(c) \leq 4^{k_l-1}\delta+3.4^{k_l-1}\delta+4^{k_{l-1}-1}\delta \leq 4^{k_l+k_{l-1}}\delta.$$

Consider the surface with boundary $\Sig'$ obtained by slitting open $\Sig$ along the trees $A_1,\dots,A_l$. Let $C_j, \; j=1,\dots,l,$ denote
the connected component of $\partial\Sig'$ arising from $A_j$, and $V'_j$ denote the finite subset of $C_j$ corresponding to the vertices of
$A_j$. We denote by $\mathrm{V}'$ the finite subset of $\Sig'$ arising from $\mathrm{V}$, note that $V'_j=\mathrm{V}'\cap C_j$. Let $\di'$
denote the distance induced by the metric structure of $\Sig'$.\\

\noindent The path $c$ represents then a path $c'$ in $\Sig'$ joining a point $x'_l$ in $V'_l$ to a point $x'_{l-1}$ in $V'_{l-1}$. Since
$\leng(c')=\leng(c) \leq 4^{k_l+k_{l-1}}\delta$, we deduces that

$$\di'(x'_{l},x'_{l-1}) \leq 4^{k_l+k_{l-1}}\delta.$$

Let $c'_0$  be a path realizing the distance from $x'_{l-1}$ to $x'_l$ in $\Sig'$, then $c'_0$ is a union of geodesic segments with endpoints in
$\mathrm{V}'$, and $\leng(c'_0)\leq 4^{k_l+k_{l-1}}\delta$.  Now, the path $c'_0$ corresponds to a path $c_0$ in $\Sig$, joining a vertex of
$A_l$ to a vertex of $A_{l-1}$. By construction, $c_0$ is a union of geodesic segments with endpoints in $\mathrm{V}$, each of which is either
an edge of a tree in $\{A_1,\dots,A_l\}$, or a geodesic segment which does not cross any edge of the trees in the family $\{A_1,\dots,A_l\}$. As
a consequence, the union of $c_0$ and all the trees in $\{A_1,\dots,A_l\}$ which have at least a common point with $c_0$ is a geodesic tree.
This new tree contains obviously $A_{l-1}$ and $A_l$ as subtrees, hence it contains at least $k_l+k_{l-1}+1$ edges. We denote by $A'_{l'}$ this
new tree, and  by $A'_1,\dots,A'_{l'-1}$  the remaining trees in the family $\{A_1,\dots,A_l\}$.\\

\noindent  It is a routine to verify that the family $\{A'_1,\dots, A'_{l'}\}$ also satisfies the conditions $a)$, and $b)$ of $(**)$. If the
condition $c)$ still holds, then we can restart the procedure. Since the number of singularities of $\Sig$ is finite, the procedure can be
repeated until we get either

\begin{itemize}
\item[.] a single geodesic tree $A$ verifying the property $(*)$ or,

\item[.] a family $\{\tilde{A}_1,\dots,\tilde{A}_{\tilde{l}}\}$ of disjoint geodesic trees, verifying $a)$, and $b)$ of the condition $(**)$,
and in addition, we have:

$$\di(\tilde{A}_{\tilde{l}},\tilde{A}_1\sqcup \dots \sqcup \tilde{A}_{\tilde{l}-1}) \geq 3.4^{k_{\tilde{l}}-1}\delta,$$

\noindent where $k_{\tilde{l}}$ is the number of edges of $\tilde{A}_{\tilde{l}}$.

\end{itemize}

Now, let us show that the tree joining procedure, and the vertex adding procedure in Lemma \ref{FVSSlm2} will allow us to construct a good
forest in $\Sig$. First, by Lemma \ref{FVSSlm2}, we know that, there exists a good tree $A_1$. We will proceed by induction. Assume that we
already have a family $\{A_1,\dots,A_l\}$ of disjoint good trees. If the union of the vertex sets of $A_1,\dots,A_l$ is $\mathrm{V}$, or all the
remaining singularities have cone angle in $2\pi\N$, then we are done. Otherwise, we can start a vertex adding procedure with a singular point
which is not a vertex of the family $\{A_1,\dots,A_l\}$. \\

\noindent The vertex adding procedure can be carried out until we get a new good tree $A_{l+1}$ disjoint from $A_1\sqcup\dots\sqcup A_l$, or
until we get a geodesic tree $A$ such that

\begin{itemize}

\item[.] $A$ satisfies the condition $(*)$,

\item[.] the segment realizing the distance $\di(\ver(A),\mathrm{V}\setminus\ver(A))$ intersects a tree in the family $\{A_1,\dots,A_l\}$.\\

\end{itemize}

\noindent In the latter case, we see that the family $\{A_1,\dots,A_l,A\}$ satisfies the condition $(**)$, therefore we can start  the tree
joining procedure. When this procedures terminates, we get a family of disjoint geodesic trees $\{\tilde{A}_1,\dots,\tilde{A}_{\tilde{l}}\}$, it
may happen that $\tilde{l}=1$, where $\tilde{A}_1,\dots,\tilde{A}_{\tilde{l}-1}$ are good, $\tilde{A}_{\tilde{l}}$ verifies the condition $(*)$,
and we can carry out the vertex adding procedure on $\tilde{A}_{\tilde{l}}$. Since the number of singularities of $\Sig$ is finite, this
algorithm must terminate, and we obtain a good forest for $\Sig$.\carre

\section{Proof of Proposition \ref{FVSSpr2}}\label{FVSSpr2Sect}

Let $I_0$ be a subset of $\{1,\dots,n\}$ such that $\delta^+_{I_0}(\Sig)=\delta(\Sig)=\delta$. Let $s$ be a geodesic segment joining a point
$x_{i_0}$ with  $ i_0\in I_0$ and a point $x_{i_1}$ with $i_1 \not\in I_0 $ such that $\leng(s)=\delta$. Let $p$ denote the midpoint of $s$. As
usual, we denote by $\di$ the distance induced by the flat metric of $\Sig$. First, we have

\begin{lemma}\label{FVSSlm4}
$B(p,\delta/2)=\{x\in \Sig: \di(p,x)<\delta/2\}$ does not contain any singular point of $\Sig$.\\
\end{lemma}

\dem Suppose on the contrary that a singular point  $x_k$, with $k\not\in\{i_0,i_1\}$, is contained in $B(p,\delta/2)$, then we have
$\di(x_{i_0},x_k) < \delta$, and $\di(x_{i_1},x_k)< \delta$, but this would imply that $\delta_{I_0}(\Sig)<\delta$, and we have a
contradiction.\carre

Let $D(\delta/2)$ denote the open disk with center $(0,0)$ and radius $\delta/2$ in the Euclidean plane $\E^2=\R^2$. Let $f$ be the isometric
immersion from $D(\delta/2)$ to $\Sig$, which maps the horizontal diameter of $D(\delta/2)$ to the segment $s$, and the origin $(0,0)$ to the
point $p$. The immersion $f$ exists because the smallest distance from $p$ to a singular point of $\Sig$ is $\delta/2$.\\

\noindent Let $\epsilon$ be the maximal value such that the restriction of $f$ on the disk $D(\epsilon\delta)$ with center $(0,0)$ and radius
$\epsilon\delta$ is an embedding. If $\epsilon\geq 1/4$ then there is an embedded Euclidean disk of radius $\delta/4$ in $\Sig$, which means
that $\mathbf{Area}(\Sig) \geq (\pi\delta^2)/{16}$. In what follows, we will suppose that $\epsilon <1/4$, consequently, the set $f^{-1}(\{p\})$
contains points other than $(0,0)$. Let $p_1$ be the point in $f^{-1}(\{p\})\setm\{(0,0)\}$ closest to $(0,0)$.\\

For any subset $I$ of $\{1,\dots,n\}$, we denote by $\alpha_I$ the sum $\sum_{i\in I}\alpha_i$, and $\| \alpha_I\|$ the distance from $\alpha_I$
to the set $\pi\Z$ in $\R$. Set

$$\alpha_0=\min\{\|\alpha_I\|: \; I \subset \{1,\dots,n\}, \|\alpha_I\| \neq 0\}.$$

\noindent Choose a number $\epsilon_0$ such that

$$\epsilon_0 < \min\{1/6, \sin(\alpha_0)/4\}.$$

\noindent We will prove that there exists an embedded disk of radius $\epsilon_0\delta$ in $\Sig$, which is enough to prove the proposition.\\

Let $d_0$ denote the horizontal diameter of $D(\delta/2)$, and $d_1$ denote the lift of $s$ passing through $p_1$. Let $c_1$ denote the segment
joining $(0,0)$ to $p_1$ in $D(\delta/2)$, and $c$ denote the image of $c_1$ under $f$, $c$ is then a geodesic loop in $\Sig$ with base point
$p$. Let $\theta$ be angle between $d_0$ and $d_1$, by this we mean the angle in $[0;\pi/2]$ between the two lines supporting $d_0$ and $d_1$.
First, let us prove

\begin{lemma}\label{FVSSlm5}
We have either $\theta=0$, or $\epsilon >\epsilon_0$.\\
\end{lemma}

\dem Remark that  $\theta$ equals the rotation angle of the holonomy of $c$, which is the sum of some angles in $\{\alpha_1,\dots,\alpha_n\}$
modulo $\pi$. Suppose that $\theta\neq 0$, then, by the definition of $\alpha_0$, we have $\theta\geq \alpha_0$.\\

\noindent  If $\epsilon < \epsilon_0$, then the distance from $(0,0)$ to $d_1$ is less than $2\epsilon_0\delta<\sin(\alpha_0)\delta/2$. Together
with the fact that $\theta \geq\alpha_0$, this implies that $d_1$ intersects $d_0$, in other words, the segment $s$ has self-intersection, which
is impossible. Therefore, we can conclude that either $\theta=0$, or $\epsilon>\epsilon_0$.\carre

\vspace{0.5cm}

If $\epsilon > \epsilon_0$, then we are done. Therefore, we only have to consider the case $\theta=0$, and we have\\

\begin{lemma}\label{FVSSlm6}
If $\theta=0$, then the rotation angle of the holonomy of $c$ is $0$ modulo $2\pi$.\\
\end{lemma}

\dem If it is not the case, then this angle equals $\pi$ modulo $2\pi$, and hence, the holonomy of $c$ is the composition of a rotation of angle
$\pi$ and a translation which maps $(0,0)$ to $p_1$. Such a transformation must fix the midpoint $q_1$ of the segment joining $(0,0)$ to $p_1$.
It follows that $q_1$ is mapped by $f$ into a singular point of $\Sig$, which is impossible because $q_1$ is contained in the disk
$D(\delta/2)$. \carre

From Lemma \ref{FVSSlm6}, we deduce that the image of $D(\delta)$ under $f$ contains a cylinder $C$ with length $(1-2\epsilon)\delta$ and width
bounded by $2\epsilon\delta$. Remark that $c$ is then a closed geodesic in $C$ which cuts $\Sig$ into two flat surfaces with geodesic boundary,
each of which is homeomorphic to a topological closed disk. We denote by $\Sig_0$ the flat disk that contains $x_{i_0}$.\\

\begin{lemma}\label{FVSSlm7}
For any  $i$ in $I_0$,  $x_i$ is contained in $\Sig_0$.\\
\end{lemma}

\dem  Recall that by the definition of $\delta$, we have

$$\diam\{x_i, \; i\in I_0 \} < \delta/3,$$

\noindent  which implies that $\di(x_{i_0},x_i)<\delta/3$, for any $i$ in $I_0$. If there exists $i\in I_0$ such that $x_i \not\in \Sig_0$, then
the path realizing the distance $\di(x_{i_0},x_i)$ must intersect the closed geodesic $c$, therefore it crosses $C$. Consequently,

$$ \di(x_{i_0},x_i) \geq (1-2\epsilon)\delta > 2/3\delta,$$

\noindent which is impossible.\carre

\vspace{0.5cm}

The rotation angle of the holonomy of $c$ equals the sum of all cone angles at singular points in $\Sig_0$ modulo $2\pi$. By assumption, we know
that $\alpha_{I_0}\not\in 2\pi\Z$, it means that $\Sig_0$ contains singular points which do not belong to $\{x_i, \; i\in I_0\}$. Note that we
have

$$\min\{\di(x_i,x_j\}, \; i\in I_0,\; j\not\in I_0,\; x_j\in \Sig_0\} \geq \delta_{I_0}(\Sig)=\delta.$$

\noindent Since $\Sig_0$ is a flat surface with geodesic boundary which contains no singularities on the boundary, we can restrict ourselves
into $\Sig_0$ and restart the whole procedure.  This procedure can be continued as long as  the rotation angle of the loop $c$ is zero.\\

\noindent Since we only have finitely many singular points in $\Sig$, the procedure must stop, and we get a point in $\Sig$ whose injectivity
radius is at least $\epsilon_0\delta$. Proposition \ref{FVSSpr2} is then proved. \carre

\end{appendices}


\end{document}